\documentclass[11pt]{scrartcl}
\usepackage[T1]{fontenc}
\usepackage[utf8]{inputenc}	

\usepackage{amsfonts,amssymb,amsmath}
\usepackage{multimedia}
\usepackage{float}
\usepackage{placeins}
\usepackage{comment} 
\usepackage{url}

\usepackage{enumitem}
\usepackage{mathabx}

\usepackage{tikz} 
\usepackage{titletoc}
\usepackage{subfig}
\usepackage[subfigure]{tocloft}

\definecolor{tumo}{RGB}{227,114,34}
\definecolor{tumg}{RGB}{218,215,203}
\definecolor{tumb}{RGB}{0,101,189} 
\definecolor{tumd}{RGB}{0,82,147} 
\definecolor{tumblues2}{RGB}{0,51,89}
\definecolor{headercolor}{RGB}{0,0,0} 
\definecolor{tumgray2}{gray}{0.5}

\usepackage{tabularx}

\newlength\origtabcolsep
\origtabcolsep=\tabcolsep
\newcolumntype{e}{>{\hsize=0pt}X}
\newcolumntype{x}{>{\hskip\origtabcolsep}X<{\hskip\origtabcolsep}}
\newcolumntype{C}{>{\hskip\origtabcolsep}c<{\hskip\origtabcolsep}}
\newcolumntype{R}{>{\hskip\origtabcolsep}r<{\hskip\origtabcolsep}}

\usepackage{longtable}
\usepackage{bm} 
\usepackage{mathtools} 
\usepackage{bbm} 
\usepackage{graphicx} 
\usepackage{multirow}

\usepackage{colortbl}

\makeatletter
\DeclareOldFontCommand{\rm}{\normalfont\rmfamily}{\mathrm}
\DeclareOldFontCommand{\sf}{\normalfont\sffamily}{\mathsf}
\DeclareOldFontCommand{\tt}{\normalfont\ttfamily}{\mathtt}
\DeclareOldFontCommand{\bf}{\normalfont\bfseries}{\mathbf}
\DeclareOldFontCommand{\it}{\normalfont\itshape}{\mathit}
\DeclareOldFontCommand{\sl}{\normalfont\slshape}{\@nomath\sl}
\DeclareOldFontCommand{\sc}{\normalfont\scshape}{\@nomath\sc}
\makeatother

\newtheorem{remark}{Remark}

\DeclareMathAlphabet{\pazocali}{OMS}{zplm}{m}{n}

\newcommand{\mathcaa}{\pazocali}

\newcommand{\myquad}{\hspace*{4ex}}

\newcommand{\qfct}{{\mathbf q}}

\newcommand{\mmtau}{{  \sdom{ \boldsymbol \tau}}}

\newcommand{\adom}{}
\newcommand{\aledom}{\hat }
\newcommand{\sdom}{\tilde}
\newcommand{\pdom}{\check}
\newcommand{\widealedom}{\widehat}
\newcommand{\widesdom}{\widetilde}

\newcommand{\mId}{\mathbf{I}}
\newcommand{\mdl}{(\!(}
\newcommand{\mdr}{)\!)}

\newcommand{\mlnorm}{\|}
\newcommand{\mrnorm}{\|}

\newcommand{\mx}{{\sf{x}}}
\newcommand{\my}{{\sf{y}}}
\newcommand{\mz}{{ \sf{z}}}

\newcommand{\mt}{{\sf{t}}}
\newcommand{\ms}{{\sf{s}}}

\newcommand{\mevalt}{{(\cdot,\mt)}}
\newcommand{\mevals}{{(\cdot,\ms)}}

\newcommand{\mevalo}{{(\cdot,\sf {0})}}

\newcommand{\mnormalf}{{\mathbf n_f}}
\newcommand{\mnormals}{{\mathbf n_s}}
\newcommand{\manf}{\adom{\mathbf n}_f}

\newcommand{\malenf}{\aledom{\mathbf n}_f}
\newcommand{\malens}{\aledom{\mathbf n}_s}

\newcommand{\mQf}{{Q_f^T}}
\newcommand{\mQs}{{Q_s^T}}
\newcommand{\mGf}{\Sigma_f^T}
\newcommand{\mGs}{ \Sigma_s^T}
\newcommand{\mGi}{\Sigma_i^T}
\newcommand{\mhQf}{{\aledom Q_f^T}}
\newcommand{\mhQs}{{\aledom Q_s^T}}

\newcommand{\mhGfN}{{\aledom \Sigma^T_{fN}}}
\newcommand{\mhGsN}{{\aledom \Sigma^T_{sN}}}

\newcommand{\mDy}{D_\my\,}
\newcommand{\mDz}{D_\mz\,}

\newcommand{\msigmaf}{\boldsymbol \sigma_f}
\newcommand{\msigmas}{\boldsymbol \sigma_s}

\newcommand{\mTadh}{{  \sdom{\bm{\mathcaa T}} _{ad, \sf h}}}

\newcommand{\mX}{{\boldsymbol \chi}}

\newcommand{\mXs}{{\boldsymbol \chi_s}}

\newcommand{\p}{{p}}
\newcommand{\vel}{{\mathbf v}} 
\newcommand{\vela}{{\mathbf u}} 

\newcommand{\dis}{{\mathbf w}} 
\newcommand{\disa}{{\mathbf z}} 

\newcommand{\mv}{\vel}

\newcommand{\mrhof}{\rho_f}
\newcommand{\mrhos}{\rho_s}

\newcommand{\mpsiv}{{\boldsymbol \psi}^v}
\newcommand{\mpsip}{{\psi}^p}
\newcommand{\metav}{{\boldsymbol \eta}^v}
\newcommand{\metap}{{ \eta}^p}
\newcommand{\mpsiw}{{\boldsymbol \psi}^z}
\newcommand{\mpsiu}{{\boldsymbol \psi}^w}

\newcommand{\mbetav}{\overline{\boldsymbol \eta}^v}
\newcommand{\mbetap}{\overline{ \eta}^p}
\newcommand{\mbpsiv}{\overline{\boldsymbol \psi}^v}
\newcommand{\mbpsip}{\overline{ \psi}^p}

\newcommand{\maledomJXs}{\widehat{\mathbf F}_{\mXs}}
\newcommand{\maledomJX}{{\widehat{\mathbf F}_\mX}}

\newcommand{\msdomJXt}{{\widetilde{\mathbf F}_{\mX}}}

\newcommand{\maledomJ}{{\aledom{J}_\mX}}
\newcommand{\mFX}{\msdomJXt}
\newcommand{\mFXa}{\maledomJX}
\newcommand{\msdomJ}{\sdom{J}_\mX}

\newcommand{\mU}{\mathbf W  (\mv )}
\newcommand{\mTad}{{  \sdom{\bm{\mathcaa T}} _{ad}}}

\newcommand{\maledomEs}{{\widealedom{\mathbf E}_{\mXs}}}
\newcommand{\msdomE}{{\widesdom{\mathbf E}_\mX}}
\newcommand{\mmut}{ \sdom {\mathbf u }_\tau }

\newcommand{\mUt}{{\sdom{\mathbf U}_{ad}}}
\newcommand{\mtuht}{\sdom \dis} 

\newcommand{\mtJht}{\msdomJ}
\newcommand{\mtJtt}{\mathrm{det}\mtFtt}
\newcommand{\mtFtt}{ ({D_\mz \mmtau} )}
\newcommand{\mtFti}{ (D_\mz \mmtau )^{-1}}
\newcommand{\mtFi}{ ({D_\mz \mmtau} )^{-\top}}
\newcommand{\mtFht}{\mFX}

\title{Numerical Methods for Shape Optimal Design of Fluid-Structure Interaction Problems}
\author{Johannes Haubner \thanks{Department of Mathematics and Scientific Computing, University of Graz, Heinrichstr. 36, 8010 Graz, Austria \mbox{(johannes.haubner@uni-graz.at)}}
\and Michael Ulbrich \thanks{Department of Mathematics, Technical
	University of Munich, Boltzmannstr.~3, 85748 Garching b. M\"unchen,
	Germany \mbox{(mulbrich@ma.tum.de)}}
}

\pgfdeclarelayer{bg1}
\pgfdeclarelayer{bg2}
\pgfdeclarelayer{bg3}
\pgfdeclarelayer{bg4}
\pgfsetlayers{bg1, bg2, bg3,bg4,  main}

\begin{document}
	
\maketitle

\begin{abstract}
We consider the method of mappings for performing shape optimization for unsteady fluid-structure interaction (FSI) problems.
In this work, we focus on the numerical implementation. 
We model the optimization problem such that it takes several theoretical
results into account, such as regularity requirements on the transformations and a differential geometrical point of view on 
the manifold of shapes. 
Moreover, we discretize the problem such that we can compute exact discrete gradients. 
This allows for the use of general purpose optimization solvers.
We focus on an FSI benchmark problem to validate our numerical implementation. 
The method is used to optimize parts of the outer boundary and the interface.
The numerical simulations build on FEniCS,
dolfin-adjoint and IPOPT. Moreover, as an additional theoretical result, 
we show that for a linear special 
case the adjoint attains the same structure as the forward problem but reverses the 
temporal flow of information.
\end{abstract}

\section{Introduction}
Fluid-structure interaction (FSI) is a particularly important subclass of multi-physics problems that arise frequently in applications such as wind turbines, bridges, naval architecture or biomedical applications, cf., e.g., 
\cite{DoHa}.
	Performing shape optimization on this class of problems is also a vivid area of 
	research. In \cite{HUU19} this task is tackled from a theoretical point of view, 
	whereas
	\cite{HeDuBl, HRHD18, HoStGaIsWuBl, LoPaQuRo, LMB14, LuMoJa, MaAlRe, ML13, SNGWB17} focus on numerical approaches. In this paper, we present our numerical realization of shape optimization for an unsteady nonlinear FSI problem.

\subsection{Method of mappings}

The different methods for performing shape optimization are closely related to the metric structure 
that is imposed on the set of admissible shapes. Besides the consideration of characteristic functions (which motivates, e.g., phase field approaches, cf., e.g., \cite{GHHKL16}, or \cite{HNU23}) or distance functions \cite{DeZo11}, a metric can be defined via transformations \cite{MuSi}. Similarly to the FSI problem, the Lagrangian or Eulerian perspective can be chosen to work with transformations. The latter leads to the notion of shape derivatives and to level set methods. The Lagrangian perspective is called method of mappings or perturbation of the identity 
\cite{MuSi, BeFe, BaLiUl, GuMa, KeUl}
. It allows for a rigorous theoretical framework for shape optimization and was already applied in fluid
	mechanics, see e.g. \cite{BaLiUl}. Instead of optimizing over a set of admissible shapes, 
a parametrization of the shape by
a bi-Lipschitz homeomorphism
$\mmtau: \mathbb R^d \to \mathbb R^d$ via $ \aledom \Omega
= \mmtau(   \sdom \Omega )$ is used, where $  \sdom \Omega \subset \mathbb R^d$ is a
shape reference domain the coordinates of which we will denote by $\mz$. The optimization is then 
performed over a set of admissible transformations. For more details we refer to the cited literature 
and Section \ref{sec::fsisop}.
In this work, we apply the method of mappings to solve a shape optimization problem for unsteady nonlinear FSI.

\subsubsection{Theoretical challenges}
A main theoretical challenge is ensuring existence of solutions of shape optimization problems, which, in general, requires a regularization term and can be done with available optimization theory since the problem is formulated on a fixed reference domain. If a partial differential equation (PDE) constraint is involved, one crucial and technical step is proving differentiability of the state with respect to shape transformations \cite{FLUU16, HUU19}. In this work, we do not focus on the theoretical aspects of shape optimization of fluid-structure interaction but concentrate on the algorithmic realization.

\subsubsection{Algorithmic realization}
\label{secAlgReal}
Treating, for the moment, the computation of the surface or volumetric shape functional and its derivative as black-box, our proposed algorithm takes the following form:
\begin{enumerate}
\item Choose a design boundary $\tilde \Gamma_d$ of the initial domain $\tilde \Omega$. 
\item Choose a scalar valued quantity $\tilde d: \tilde \Gamma_d \to \mathbb R$ as optimization variable. 
\item Apply a boundary operator $S_\Gamma$ to obtain a vector valued deformation field on the boundary. 
\item Apply a suitable extension operator $S_\Omega$ to $S_\Gamma(\tilde d)$, such that a deformation field for $\tilde \Omega$ is given by $S_\Omega(S_\Gamma(\tilde d))$. For abbreviation, we introduce $\mathrm{B} := S_\Omega \circ S_\Gamma$.
\item Compute the value of the shape functional and its derivative at $(\mathrm{id} + \mathrm{B}(\tilde d))(\tilde \Omega)$. 
\item Use the chain rule to obtain the derivative with respect to $\tilde d$. 
\item Use this within an appropriately chosen derivative-based optimization algorithm. If second order derivatives are available, also second order optimization methods can be applied.
\end{enumerate}

\subsubsection{Requirements on software framework}

The above algorithm requires that the software framework allows for the solution of PDEs on (sub)meshes and (submeshes of) their surface meshes. Moreover, we assume that it is possible to restrict functions on the mesh to functions on submeshes and (submeshes of) their surface meshes and vice versa extend functions that live on parts of the surface appropriately (e.g. by zero) to functions on the mesh. Our implementation builds on FEniCS \cite{ FenicsBook}.

\subsubsection{Related work}

One common approach to deal with transformation based shape optimization is solving a sequence of auxiliary problems with linearized objective function based on a shape derivative to iteratively update the domain, e.g., \cite{radtke2023parameter, deckelnick2022novel, schulz2016efficient} and references therein. In this work, we solve the shape optimization problem on a fixed domain ($\tilde \Omega$) without linearization of the objective function.
The decision for a scalar valued optimization variable is motivated by considerations on one-to-one correspondence between shapes and controls, see, e.g., \cite{HSU20}.

The choice of the operator $\mathrm B$ is part of the modeling of the shape optimization problem, see, e.g., \cite{HSU20} and references therein. 
In some approaches, Steklov-Poincar\'e operators are used, e.g., \cite{schulz2016efficient}, which can be incorporated into our framework via $S_\Gamma$. Other approaches involve remeshing or a multi-mesh approach \cite{dokken2019shape}, which can also be incorporated in our approach, e.g. via $S_\Omega$ if only the interior of the domain is remeshed. In this paper, we want to work with exact (discrete) derivatives during the optimization process. Therefore, we do not consider remeshing routines. Moreover, in order to have full control over the shape functional and its derivative, we do not consider the computation of the derivative as black box.

\subsection{Fluid-Structure Interaction Problem}
We consider FSI problems for which the fluid is modeled by the unsteady incompressible Navier-Stokes equations. These are canonically formulated in the Eulerian framework, i.e., on a time-dependent physical domain $\pdom \Omega_f (\mt) \subset \mathbb R^d$ for $\mt \in I:= (0,T)$, $T>0$. We divide the fluid boundary $\partial \pdom \Omega_f (\mt) = \pdom \Gamma_{fD}(\mt) \cup \pdom \Gamma_{fN}(\mt) \cup \pdom \Gamma_i (\mt)$ into three disjoint parts on which Dirichlet (on $\pdom \Gamma_{fD}(\mt)$) or Neumann (on $\pdom \Gamma_{fN}(\mt)$) boundary conditions, or coupling conditions (on the interface $\pdom \Gamma_i(\mt)$) are imposed. The corresponding space-time cylinders are denoted by \begin{align*} &\pdom Q_f^T := \bigcup_{\mt \in I} \pdom \Omega_f(\mt)\times \lbrace \mt \rbrace, \quad \pdom \Sigma_{fD}^T := \bigcup_{\mt \in I} \pdom \Gamma_{fD}(\mt)\times \lbrace \mt \rbrace, \quad \\ &\pdom \Sigma_{fN}^T := \bigcup_{\mt \in I} \pdom \Gamma_{fN}(\mt)\times \lbrace \mt \rbrace, \quad \pdom \Sigma_{i}^T := \bigcup_{\mt \in I} \pdom \Gamma_{i}(\mt)\times \lbrace \mt \rbrace.
\end{align*}
The differential equations are given by
\begin{align*}
\hspace{12ex}\rho_f \partial_\mt \pdom \mv_f  + \rho_f (\pdom \mv_f \cdot \nabla_\mx ) \pdom \mv_f -  \mathrm{div}_\mx ( \sigma_{f,\mx}(\pdom \mv_f ,  \pdom \p_f)) &= \rho_f \pdom{\mathbf f}_f && \text{in } \pdom Q_f^T,\hspace*{12ex} \\
\mathrm{div}_\mx (\pdom \mv_f) &=0&& \text{in } \pdom Q_f^T , \\
\pdom \mv_f &= \pdom \mv_{fD}&& \text{on } \pdom \Sigma_{fD}^T, \\
\sigma_{f,\mx} ( \pdom \mv_f , \pdom p_f ) \pdom{ \mathbf n}_f  &= \pdom{ \mathbf g}_f  && \text{on } \pdom \Sigma_{fN}^T,
\end{align*}
with the initial condition
\[
\pdom \mv_{f} \mevalo = \pdom \mv_{0f} \myquad \text{on } \pdom\Omega_f(\sf{0}),
\]
where $\pdom \mv_f$ denotes the fluid velocity, $\pdom p_f$ the pressure and $\pdom{ \mathbf n}_f$ the outer unit normal vector. $\pdom{\mathbf f}_f$, $\pdom \mv_{fD}$, $\pdom{ \mathbf g}_f$ and $\pdom \mv_{0f}$ are right-hand side, boundary and initial values. The fluid stress tensor is defined by 
\begin{align*}
\sigma_{f,\mx} (\pdom \mv_f, \pdom p_f ) =  \rho_f \nu_f  (D_\mx \pdom \mv_f + {D_\mx \pdom \mv_f}^T  )- \pdom p_f \mId,
\end{align*}
with unit matrix $\mId \in \mathbb{R}^{d\times d}$ and $D_\mx \pdom\mv := (\partial_{\mx_j} \pdom\mv_i)_{i,j}$ denoting the Jacobian of $\pdom\mv$. The parameters $\rho_f$ and $\nu_f$ denote the fluid density and viscosity, respectively. The structure equations, however, are formulated in the Lagrangian framework, i.e., on a fixed reference domain $\aledom \Omega_s$ with disjoint Dirichlet, Neumann, and interface boundary parts $\aledom \Gamma_{sD}$, $\aledom \Gamma_{sN}$, and $\aledom \Gamma_i$ such that $\partial \aledom \Omega_s = \aledom \Gamma_{sD} \cup \aledom \Gamma_{sN} \cup \aledom \Gamma_i$. The physical domain $\pdom \Omega_s (\mt)$ for any $\mt \in I$ is obtained by the transformation $\aledom \mX_s \mevalt: \aledom{ \Omega}_s \to \pdom \Omega_s (\mt )$, $\aledom{\mX }_s ( \my, t )= {\my} + \aledom{\dis}_s ({\my},\mt )$, where the deformation $\aledom{ \dis}_s$ solves the hyperbolic equations
\begin{align*}
\rho_s  \partial_{\mt\mt} \aledom{\dis}_{s} - \mathrm{div}_\my ( \maledomJXs {\Sigma}_{s,\my} ( \aledom{ \dis}_s )  )= \rho_s \aledom{\mathbf f}_s & \quad \text{in }\aledom Q_s^T := \aledom\Omega_s\times I, \\
\aledom{\dis}_{s} = \aledom{\dis}_{sD}& \quad \text{on } \aledom \Sigma_{sD}^T := \aledom\Gamma_{sD} \times I,\\
\maledomJXs{\Sigma}_{s,\my} ( \aledom{\dis}_s  ) \aledom{\mathbf n}_s = \aledom{\mathbf g}_s&\quad \text{on } \aledom \Sigma_{sN}^T := \aledom\Gamma_{sN} \times I,\\
\aledom{\dis}_{s} \mevalo = \aledom{\dis}_{0s} & \quad \text{on } \aledom \Omega_s, \\ \partial_\mt  \aledom{\dis}_{s} \mevalo = \aledom{\dis}_{1} & \quad \text{on } \aledom \Omega_s,
\end{align*}
and we define $\maledomJXs := D_\my \aledom{\mX}_s$.
Here, $\rho_s$ denotes the structure density and $\aledom{\mathbf f}_s$, $\aledom{\dis}_{sD}$, $\aledom{ \mathbf g}_s$, $\aledom{\dis}_{0s}$ and $\aledom{\dis}_{1}$ denote right hand side, boundary and intial values. For nonlinear Saint Venant Kirchhoff type material the stress tensor ${ \Sigma}_{s,\my} ( \aledom{\dis}_s  )$ is given by 
\begin{align*}
{{\Sigma}}_{s,\my} ( \aledom{\dis}_s  ) := \lambda_s \mathrm{tr} ({\maledomEs}  )\mathbf I + 2 \mu_s {\maledomEs}, 
\end{align*}
with ${\maledomEs}  := \frac{1}{2} ({\maledomJXs}^\top {\maledomJXs}  - \mId  )$.
The first challenge for considering the coupled problem arises from the fact that the above canonical models for the fluid and structure equations are formulated in different frameworks for which reason a precise formulation of the coupling conditions requires further effort.

For FSI simulations, partitioned as well as monolithic approaches have been proposed. Partitioned methods solve the corresponding models seperately and typically apply fixed point iterations to the coupling interface conditions, which can, e.g., be accelerated by Quasi-Newton techniques \cite{DeBaVi, KaIbNiMa}.  
Monolithic approaches 
\cite{DoHuPoRo, DuRaRi, FrRiWi1, GhLi, HeHaBo, Wi}
, such as arbitrary Lagrangian-Eulerian (ALE) \cite{DoHuPoRo, DuRaRi, HeHaBo} and fully Eulerian methods 
\cite{DuRaRi, FrRiWi1, Wi}
, use the same reference frame for fluid and solid.  While fully Eulerian approaches use the spatial reference frame, the ALE framework is obtained by introducing an arbitrary but fixed reference domain $\aledom \Omega_f$ such that the fluid and solid reference domains are disjoint, i.e., $\aledom \Omega_s \cap \aledom \Omega_f = \emptyset$, and share the interface $\aledom \Gamma_i := \overline{\aledom \Omega_s} \cap \overline{\aledom \Omega_f}$ as part of their boundaries. Moreover, an extension $\aledom{\mX}(\cdot, \mt) : \aledom \Omega \to \mathbb R^d$ of the solid transformation $\aledom{\mX}_s(\cdot, \mt) $ to the whole reference domain $\aledom \Omega := \aledom \Omega_s \cup \aledom \Omega_f \cup \aledom \Gamma_i$ is introduced for any $\mt \in I$, i.e., $\aledom{\mX}(\cdot, \mt)\vert_{\aledom \Omega_s} = \aledom{\mX}_s(\cdot, \mt)$. It can, e.g., be obtained by choosing a fully Lagrangian setting or a harmonic or biharmonic extension of the solid displacement to the fluid reference domain. Transformation of the fluid equations with the help of $\aledom {\mX}$ to the fixed reference domain $\aledom \Omega_f$ and coupling the fluid and structure equations across the interface $\aledom \Gamma_i$ yields the system of equations
\begin{align}
\begin{split}
{\maledomJ} \rho_f \partial_\mt \aledom{\mv}_f + {\maledomJ} \rho_f  ( ({\maledomJX}^{-1}  (\aledom{ \mv}_f - \partial_\mt \aledom{\mX}  ) )\cdot \nabla_\my ) \aledom{\mv}_f & \\ - \mathrm{div}_\my ({\maledomJ} \aledom{\boldsymbol \sigma}_{f} {\maledomJX}^{-\top}  ) = {\maledomJ} \rho_f \aledom{\mathbf f}_f& \quad \text{in } \aledom Q_f^T := \aledom\Omega_f \times I,\\
\mathrm{div}_\my  ({\maledomJ} {\maledomJX}^{-1} \aledom{\mv}_f ) = 0 & \quad \text{in } \aledom Q_f^T, \\
\rho_s \partial_{\mt} \aledom{\mv}_s - \mathrm{div}_\my  ({\maledomJ} \aledom{\boldsymbol \sigma}_{s} {\maledomJX}^{-\top}   ) =  \rho_s \aledom{\mathbf f}_s& \quad \text{in }\aledom Q_s^T, \\
\rho_s\partial_{\mt} \aledom{\dis}_s -  \rho_s\aledom{\mv}_s = 0& \quad \text{in }\aledom Q_s^T,
\end{split}
\label{coupledsystem}
\end{align}
with initial conditions
\begin{align*}
\aledom{\mv}_f \mevalo = \aledom{\mv}_{0f} \quad\text{on } \aledom \Omega_f,\quad\quad \aledom{\dis}_s \mevalo = \aledom{\dis}_{0s} \quad \text{on } \aledom \Omega_s,\quad\quad
\aledom{\mv}_s \mevalo = \aledom{\dis}_{1} \quad \text{on } \aledom \Omega_s,
\end{align*}
boundary conditions
\begin{align*}
\aledom{\mv}_f = \aledom{\mv}_{fD}\quad \text{on } \aledom \Sigma_{fD}^T, \quad\quad \aledom{J}  \aledom{\boldsymbol \sigma}_{f} {\maledomJX}^{-\top} \aledom{\mathbf{n}}_f = \aledom{\mathbf g}_f\quad \text{on } \aledom \Sigma_{fN}^T,\\ \aledom{\dis}_{s} = \aledom{\dis}_{sD} \quad \text{on } \aledom \Sigma_{sD}^T,\quad\quad
{\maledomJ} \aledom{\boldsymbol \sigma}_{s} {\maledomJX}^{-\top} \aledom{\mathbf n}_s = \aledom{\mathbf g}_s\quad \text{on } \aledom \Sigma_{sN}^T,
\end{align*}
where $\aledom \Sigma_{fD}^T := \aledom{\Gamma}_{fD} \times I$ and $\aledom \Sigma_{fN}^T := \aledom\Gamma_{fN} \times I$,
and additional coupling conditions
\begin{align*}
 \aledom{ \mv}_s =\aledom{ \mv}_f& \quad \text{on } \aledom \Sigma_i^T := \aledom\Gamma_i \times I,\\
-{\maledomJ} \aledom{\boldsymbol \sigma}_{f} {\maledomJX}^{-\top} \aledom{\mathbf n}_f = {\maledomJ} \aledom{\boldsymbol \sigma}_s {\maledomJX}^{-\top} \aledom{\mathbf n}_s&\quad \text{on } \aledom \Sigma_i^T,
\end{align*}
where the transformed fluid stress tensor is given by
\begin{align*}
\aledom{ \boldsymbol \sigma}_{f} := \rho_f \nu_f  ( D_\my \aledom{\mv}_f\maledomJX^{-1} + \maledomJX^{-\top} D_\my \aledom{ \mv}_f^\top  ) - \aledom{p}_f \mathbf \mId .
\end{align*}
Here, $\aledom {\mathbf f }_f := \pdom {\mathbf f}_f \circ \aledom \mX$ and $\aledom {\mathbf g}_f$, $\aledom \mv_{fD}$, as well as, $\aledom \mv_{0f}$ are defined analogously. Moreover, $\aledom{\mv}_f = \pdom \mv_f \circ \aledom{ \mX}$, $\aledom p_f = \pdom p_f \circ \aledom \mX$, $\aledom{\boldsymbol \sigma}_{f} =  \sigma_{f,\mx}(\pdom \mv_f, \pdom p_f) \circ \aledom{\mX}$, $\aledom{\boldsymbol \sigma  }_{s} := {\maledomJ}^{-1}\maledomJX { \Sigma}_{s,\my}(\aledom \dis_s) {\maledomJX}^{\top}$, where ${\maledomJX} := D_\my \aledom{\mX}$ and ${\maledomJ} := \mathrm{det} ( {\maledomJX} )$, as well as, $\aledom{ \mv}_s = \partial_\mt \aledom{ \dis}_s$.
In the following, we will denote coordinates on the physical domain $\pdom \Omega$ by $\mx$ and on the reference domain $\aledom \Omega$ by $\my$ and the subscripts of the nabla-operators indicate on which variables they act on.

\begin{remark}
    For a function $\pdom \qfct: ~ \pdom \Omega \times I \to \mathbb R^n$, $n \geq 1$, $\pdom \qfct \circ \aledom \mX$ denotes the function that maps $(\my, \mt) \mapsto \pdom \qfct(\aledom \mX(\my, \mt), \mt)$.
\end{remark}

\subsection{Outline}

In Section \ref{sec::fsimod} we present the FSI model. 
Some functions in the corresponding weak formulation live on parts of the
domain and include boundary conditions at the interface. This, however, might not 
be supported by finite element toolboxes, which motivates to consider a modified 
weak formulation that only works with functions defined on the whole domain and 
boundary conditions defined on the boundary of it. Therefore, we introduce elliptic extension equations - for the pressure on the solid domain
and for the solid displacement on the fluid domain - and weight them with
parameters $\alpha_p$ and $\alpha_w$. 
Numerical tests on the FSI2 benchmark (Section \ref{sec::numres})
show the validity of this modification if 
the parameters are chosen sufficiently small.
Computing the derivative in an optimal control setting efficiently requires the solution of 
adjoint equations. For a linear special case with 
stationary interface, we show in the Appendix \ref{sec::adjcon} 
that the adjoint equations of the FSI system have a similar structure
as the forward equations if we eliminate the equation $\rho_s \partial_\mt \aledom \dis_s - \rho_s \aledom \vel_s = 0$ by defining $\aledom \dis_s$ 
via the time integral of $\aledom \vel_s$. 
Section \ref{sec::fsisop} gives a comprehensive presentation of the method of 
mappings applied to the FSI problem. We discuss the choice of admissible 
shape transformations from a continuous perspective. As objective function we choose a volume formulation for the mean fluid drag. In Section
\ref{sec::fsidis} we present the discretization of the shape optimization for
FSI via the method of mappings. We consider a continuous
approach to motivate our choice of the set of admissible transformations. 
For the discretization of the FSI equations we build on existing approaches. 
In order to demonstrate the applicability of the method of mappings for FSI problems we present numerical results (Section \ref{sec::numres}) for the FSI2 benchmark using FEniCS, dolfin-adjoint and IPOPT. IPOPT works 
with the Euclidian inner product on $\mathbb R^n$ where $n$ denotes the number of degrees of freedom. To work with the correct inner product 
inherited from the continuous perspective we apply a linear 
transformation, see Section \ref{sec::fsinr}.

\section{FSI Model for Numerical Simulations}
\label{sec::fsimod}
We consider the model (\ref{coupledsystem}), which is fully described if the ALE transformation $\aledom{\mX}: \aledom \Omega \times I \to \bigcup_{\mt \in I} (\pdom \Omega(\mt)\times \lbrace \mt \rbrace)$ is defined.
\subsection{ALE Transformation}
\label{subsec::ale}
There are several possibilities to choose the ALE transformation. These include, e.g., using a fully Lagrangian approach or extending the solid displacement to the fluid domain.
In theoretical investigations, the Lagrangian approach is often chosen 
\cite{CoSh, HUU19, IgKu, KuTu, RaVa}
, i.e., the reference domain $\aledom \Omega$ is given by the initial domain $\pdom \Omega(\sf 0)$ and the transformation is induced by the velocity field $\aledom{\mv}$. This has several advantages for the theoretical analysis. On the one hand, the contributions of the nonlinear term of the Navier-Stokes equations vanish on $\hat \Omega$. Additionally, no deformation variable on the fluid domain has to be introduced. However, it has drawbacks for numerical simulations, e.g., vortices in the flow might lead to mesh degeneration even though no solid displacement takes place. Therefore, we do not use the fully Lagrangian approach in the numerical implementation and focus on other extension techniques, which are presented below.
We construct the ALE transformation by extending the solid displacement $\aledom \dis_s$ to the fluid reference domain, denoted by $\aledom \dis_f$. We define 
\begin{align*}
\aledom \mX ( \my, \mt  ):= \my +  \aledom{\dis}_f  (\my, \sf{t} ) 
\end{align*}
for every $\my \in \aledom \Omega_f$ and ${\sf{t}} \in I$. One choice is given by the harmonic extension which, as numerical tests indicate, is prone to mesh degeneration if large mesh displacements occur \cite{Wi11}. Thus, we work with the biharmonic extension, cf., e.g., \cite{HaGr}, of the solid displacement 
to the fluid domain, which is given by 
\begin{align*}
\Delta_\my^2 \aledom{ \dis}_f = 0 &	\myquad \text{in } \aledom Q_f^T, \\
\aledom \dis_f  = 0, \quad \nabla_\my \aledom \dis_f \cdot \manf  = 0 & 	\myquad \text{on }  \aledom \Sigma_{fD}^T \cup \aledom \Sigma_{fN}^T, \\
\aledom \dis_f = \aledom \dis_s, \quad \nabla_\my \aledom \dis_f \cdot \manf  = 0 &\myquad \text{on } \aledom\Sigma_i^T.
\end{align*}
For the solution of the discretized equations, 
$H^2(\aledom \Omega_f)$-conforming finite elements are needed. 
However, these elements are not necessarily available in standard finite element toolboxes. One way to circumvent this is the weak imposition of the continuity of normal derivatives across the finite element faces using a discontinuous Galerkin approach \cite{GH09}. Another approach is the consideration of a mixed formulation of the biharmonic equation
\begin{align*}
- \Delta_\my \aledom \dis_f &= \aledom \disa_f && \text{in } \aledom Q_f^T, \\
- \Delta_\my \aledom \disa_f &= 0 && \text{in } \aledom Q_f^T, \\
\aledom \dis_f &= 0 && \text{on } \aledom \Sigma_{fD}^T \cup \aledom \Sigma_{fN}^T,\hspace*{22ex} \\
\aledom \dis_f &=  \aledom \dis_s && \text{on } \aledom \Sigma_i^T, \\
\hspace*{22ex}\nabla_\my \aledom \dis_f \cdot \manf &=  0 && \text{on } \partial \aledom \Omega_f \times I,
\end{align*}
see \cite{Wi11}.

\subsection{Strong ALE Formulation}
\label{ssec::fsimod}
A full description of the FSI equations with mixed biharmonic extension of the solid deformation to the fluid domain is given by
\begin{equation}
\begin{aligned}
{\maledomJ} \rho_f \partial_\mt \aledom{\mv}_f + {\maledomJ} \rho_f  ( (\mFXa^{-1}  (\aledom{ \mv}_f - \partial_\mt \aledom{\dis}_f  ) )\cdot \nabla_\my ) \aledom{\mv}_f\qquad& \\
- \mathrm{div}_\my ({\maledomJ} \aledom{\boldsymbol \sigma}_f \mFXa^{-\top}  ) &= {\maledomJ} \rho_f \aledom{\mathbf f}_f \quad&& \text{in } \aledom Q_f^T,\\
\mathrm{div}_\my  ({\maledomJ} \mFXa^{-1} \aledom{\mv}_f ) &= 0 && \text{in } \aledom Q_f^T, \\
\rho_s \partial_{t} \aledom{\mv}_s - \mathrm{div}_\my  ({\maledomJ} \aledom{\boldsymbol \sigma}_s \mFXa^{-\top}   ) &=  \rho_s \aledom{\mathbf f}_s && \text{in } \aledom Q_s^T, \\
\rho_s\partial_{t} \aledom{\dis}_s -  \rho_s\aledom{\mv}_s &= 0 && \text{in } \aledom Q_s^T,\\
- \Delta_\my \aledom{ \dis}_f &= \aledom \disa_f && \text{in } \aledom Q_f^T, \\
- \Delta_\my \aledom \disa_f & = 0 && \text{in } \aledom Q_f^T, 
\end{aligned}
\label{coupledsystem3}
\end{equation}
with initial and boundary conditions 
\begin{align*}
\aledom{\mv}_f \mevalo = \aledom{\mv}_{0f} \quad\text{on } \aledom \Omega_f,\quad\quad \aledom{\dis}_s \mevalo = \aledom{\dis}_{0s} \quad \text{on } \aledom \Omega_s,\quad\quad
\aledom{\mv}_s \mevalo = \aledom{\dis}_{1}  \quad \text{on } \aledom \Omega_s, \\
\aledom{\mv}_f = \aledom{\mv}_{fD} \quad \text{on }\aledom{\Sigma}_{fD}^T, \quad\quad \aledom{\dis}_{s} = \aledom{\dis}_{sD} \quad \text{on } \aledom\Sigma_{sD}^T, \\
\aledom \dis_f  = 0  	\quad \text{on } \aledom \Sigma_{fD}^T \cup \aledom \Sigma_{fN}^T, \quad\quad
\nabla_\my  \aledom \dis_f \cdot \manf = 0 \quad \text{on }   \partial  \aledom \Omega_f   \times I, \\
\aledom{J}  \aledom{\boldsymbol \sigma}_f \mFXa^{-\top} \aledom{\mathbf{n}}_f = \aledom{\mathbf g}_f \quad \text{on } \aledom\Sigma_{fN}^T, \quad\quad {\maledomJ} \aledom{\boldsymbol \sigma}_s \mFXa^{-\top} \aledom{\mathbf n}_s = \aledom{\mathbf g}_s \quad \text{on } \aledom\Sigma_{sN}^T,
\end{align*}
and additional coupling conditions
\begin{align*}
&\aledom{\mv}_s =\aledom{ \mv}_f \quad \text{on } \aledom\Sigma_i^T, \quad\quad
\aledom \dis_f =  \aledom \dis_s \quad \text{on } \aledom\Sigma_i^T \\
&-{\maledomJ} \aledom{\boldsymbol \sigma}_f \mFXa^{-\top} \aledom{\mathbf n}_f = {\maledomJ} \aledom{\boldsymbol \sigma}_s \mFXa^{-\top} \aledom{\mathbf n}_s \quad \text{on } \aledom\Sigma_i^T.
\end{align*}

\subsection{Weak ALE Formulation}
We define the function spaces 
\begin{align*}
& \aledom{\mathbf V} = \mathbf V(\aledom \Omega) \subset  \lbrace \aledom \mv \in H^1 ( \aledom \Omega)^d   ~:~  \aledom \mv\vert_{\aledom \Gamma_{fD}} = \aledom \mv_{fD}  \rbrace, \\
& \aledom {\mathbf V}_0 = \mathbf V_0 (\aledom \Omega) \subset \lbrace \aledom \mv \in H^1 ( \aledom \Omega)^d  ~:~ \aledom \mv\vert_{\aledom \Gamma_{fD}} = 0  \rbrace, \\
& \aledom {\mathbf W} = \mathbf W(\aledom \Omega) \subset  \lbrace \aledom \dis \in H^1 ( \aledom \Omega)^d  ~:~ \aledom \dis\vert_{\aledom \Gamma_{fD} \cup \aledom \Gamma_{fN} } = 0, ~ \aledom \dis\vert_{\aledom \Gamma_{sD}} = \aledom \dis_{sD}  \rbrace, \\
& \aledom{\mathbf W}_{f,s,0} := \mathbf W_{f,s,0}(\aledom \Omega) \subset  \lbrace \aledom \dis \in L^2(\aledom \Omega)^d~:~ \aledom \dis \vert_{\aledom \Omega_f} \in H^1 ( \aledom \Omega_f)^d, ~ \aledom \dis \vert_{\aledom \Omega_s} \in H^1 ( \aledom \Omega_s)^d, ~  \\
& \hspace{5,5cm} \aledom \dis\vert_{\aledom \Gamma_{fD} \cup \aledom \Gamma_{fN} } = 0, ~ \aledom \dis\vert_{\aledom \Gamma_{sD}} = 0, 
~  (\aledom \dis\vert_{\aledom \Omega_f}) \vert_{\aledom \Gamma_i} = 0  \rbrace, \\
& \aledom {\mathbf W}_0 = \mathbf W_0(\aledom \Omega)  \subset  \lbrace \aledom \dis \in H^1  ( \aledom{ \Omega} )^d  ~:~ \aledom \dis\vert_{\aledom \Gamma_{fD} \cup \aledom \Gamma_{fN}} = 0, ~ \aledom \dis\vert_{\aledom \Gamma_{sD}} = 0  \rbrace, \\
& \aledom {\mathbf Z}_f = \mathbf Z_f(\aledom \Omega) \subset H^1  (\aledom \Omega_f)^d  , \\
& \aledom{\mathbf Z}_{f,s,0} := \mathbf Z_{f,s,0}(\aledom \Omega) \subset  \lbrace \aledom \disa \in L^2(\aledom \Omega)^d~:~ \aledom \disa \vert_{\aledom \Omega_f} \in H^1 ( \aledom \Omega_f)^d, ~ \\
& \hspace{5,5cm}\aledom \disa \vert_{\aledom \Omega_s} \in H^1 ( \aledom \Omega_s)^d, ~   (\aledom \disa\vert_{\aledom \Omega_s}) \vert_{\aledom \Gamma_i} = 0  \rbrace, \\
& \aledom {\mathbf Z} = \mathbf Z(\aledom \Omega) \subset  H^1  ( \aledom \Omega)^d, \\
&\aledom P_f = P_f(\aledom \Omega) \subset  \lbrace \aledom p \in L^2 (\aledom \Omega_f  ) ~:~  \int_{\aledom \Omega_f} \aledom p d\my = 0  \rbrace, \\
& \aledom P = P(\aledom \Omega)  \subset  \lbrace \aledom p \in L^2 (\aledom \Omega  )~:~ \int_{\aledom \Omega} \aledom p d \my = 0  \rbrace.
\end{align*}
as dense subspaces. The spaces on the right hand side allow to formally write down the weak formulation. For the weak formulation to be well-defined, higher regularity 
needs to be imposed. That is the reason why we defined the function 
spaces as subspaces. We work with dense subspaces in order to stress the fact that we do impose additional boundary conditions. 

In addition, let $W_{2,q} (I, \aledom{\mathbf V}  ) :=  \lbrace \aledom \mv \in L^2  (I, \aledom {\mathbf V }) ~:~ \aledom \mv_\mt \in L^q  ( I, \aledom{\mathbf V}^*  )  \rbrace, $  where $q> 0$ and $\aledom{\mathbf V}^*$ denotes the dual space of $\aledom {\mathbf V }$. The weak formulation of \eqref{coupledsystem3} is given by: \\
Find $ ( \aledom \mv, \aledom p, \aledom \dis, \aledom \disa  ) \in W_{2,2}  (I, \aledom {\mathbf V } ) \times L^2  (I, \aledom P_f   ) \times W_{2,2} ( I,  \aledom{\mathbf W}  )\times L^2  (I, \aledom{\mathbf Z}_f  )$ such that $\aledom \mv  \mevalo= \aledom \mv_0$, $\aledom \dis  \mevalo= \aledom \dis_0$ and
\begin{align*}
&\langle \aledom A_1  ( \aledom \mv, \aledom p, \aledom \dis, \aledom \disa  ),  ( \aledom \mpsiv, \aledom \mpsip, \aledom \mpsiu, \aledom \mpsiw )  \rangle :=
( {\maledomJ} \rho_f \partial_\mt \aledom \mv, \aledom \mpsiv )_{\aledom \Omega_f} \\
&\hspace{12pt}+  ({\maledomJ} \rho_f  (  (\mFXa^{-1}  ( \aledom \mv - \partial_\mt \aledom \dis  ) )\cdot \nabla_\my  ) \aledom \mv, \aledom \mpsiv  )_{\aledom \Omega_f}  +  ( {\maledomJ} \aledom{ \boldsymbol \sigma }_f \mFXa^{-\top}, D_\my \aledom \mpsiv  )_{\aledom \Omega_f} -  (\aledom{\mathbf g}_f, \aledom \mpsiv  )_{\aledom \Gamma_{fN}} \\
&\hspace{12pt}-  ({\maledomJ} \rho_f \aledom{\mathbf f}_f, \aledom \mpsiv  )_{\aledom \Omega_f} +  ( \rho_s \partial_\mt \aledom \mv, \aledom \mpsiv  )_{\aledom \Omega_s} +  ( {\maledomJ} \aledom{\boldsymbol \sigma}_s \mFXa^{-\top}, D_\my \aledom \mpsiv  )_{\aledom \Omega_s} -  ( \aledom{\mathbf g}_s, \aledom \mpsiv  )_{\aledom \Gamma_{sN}} \\
&\hspace{12pt}-  ({\maledomJ} \rho_s \aledom{\mathbf f}_s, \aledom \mpsiv  )_{\aledom \Omega_s} +  ( \rho_s ( \partial_\mt \aledom \dis - \aledom \mv  ), \aledom \mpsiu  )_{\aledom \Omega_s} +  ( D_\my \aledom \disa, D_\my \aledom \mpsiu  )_{\aledom \Omega_f} +  ( D_\my \aledom \dis, D_\my \aledom \mpsiw  )_{\aledom \Omega_f} \\
&\hspace{12pt} -  ( \aledom \disa, \aledom \mpsiw  )_{\aledom \Omega_f} +  ( \mathrm{div}_\my  ( {\maledomJ} \mFXa^{-1} \aledom \mv ), \aledom \mpsip  )_{\aledom \Omega_f}  = 0,
\end{align*}
for all $ ( \aledom \mpsiv, \aledom \mpsip, \aledom \mpsiu, \aledom \mpsiw ) \in \aledom{\mathbf V}_0 \times \aledom P_f \times \aledom{\mathbf W}_{f,s,0}  \times \aledom{\mathbf Z}_f$ and a.e. $\mt \in I$. Since we want to work with functions defined on the whole domain, we add the equations
\begin{align*}
\aledom p &= 0 \quad \text{in } \aledom Q_s^T, \\
 - \Delta_\my \aledom \disa &= 0 \quad \text{in } \aledom Q_s^T, \quad \aledom \disa\vert_{\aledom \Omega_s} = \aledom \disa \vert_{\aledom \Omega_f}  \quad \text{on } \aledom \Sigma_i^T, \quad (\nabla_\my \aledom \disa)\vert_{\aledom \Omega_s} \cdot \malens = 0 \quad \text{on } \aledom \Sigma_{sD}^T \cup \aledom \Sigma_{sN}^T,
\end{align*}
and consider $\aledom p \in \aledom P$ and $\aledom \disa \in \aledom {\mathbf Z}$.
This leads to the weak formulation:\\
Find $ ( \aledom \mv, \aledom p, \aledom \dis, \aledom \disa  ) \in W_{2,2}  (I, \aledom {\mathbf V } ) \times L^2  (I, \aledom P   ) \times W_{2,2} ( I,  \aledom{\mathbf W}  )\times L^2  (I, \aledom{\mathbf Z}  )$ such that $\aledom \mv  \mevalo= \aledom \mv_0$, $\aledom \dis  \mevalo= \aledom \dis_0$ and
\begin{align*}
&\langle \aledom A_2  ( \aledom \mv, \aledom p, \aledom \dis, \aledom \disa  ),  ( \aledom \mpsiv, \aledom \mpsip, \aledom \mpsiu, \aledom \mpsiw )  \rangle :=
( {\maledomJ} \rho_f \partial_\mt \aledom \mv, \aledom \mpsiv )_{\aledom \Omega_f} \\
&\hspace{12pt}+  ({\maledomJ} \rho_f  (  (\mFXa^{-1}  ( \aledom \mv - \partial_\mt \aledom \dis  ) )\cdot \nabla_\my  ) \aledom \mv, \aledom \mpsiv  )_{\aledom \Omega_f}  +  ( {\maledomJ} \aledom{ \boldsymbol \sigma }_f \mFXa^{-\top}, D_\my \aledom \mpsiv  )_{\aledom \Omega_f} -  (\aledom{\mathbf g}_f, \aledom \mpsiv  )_{\aledom \Gamma_{fN}} \\
&\hspace{12pt}-  ({\maledomJ} \rho_f \aledom{\mathbf f}_f, \aledom \mpsiv  )_{\aledom \Omega_f} +  ( \rho_s \partial_\mt \aledom \mv, \aledom \mpsiv  )_{\aledom \Omega_s} +  ( {\maledomJ} \aledom{\boldsymbol \sigma}_s \mFXa^{-\top}, D_\my \aledom \mpsiv  )_{\aledom \Omega_s} -  ( \aledom{\mathbf g}_s, \aledom \mpsiv  )_{\aledom \Gamma_{sN}} \\
&\hspace{12pt}-  ({\maledomJ} \rho_s \aledom{\mathbf f}_s, \aledom \mpsiv  )_{\aledom \Omega_s} +  ( \rho_s ( \partial_\mt \aledom \dis - \aledom \mv  ), \aledom \mpsiu  )_{\aledom \Omega_s} +  ( D_\my \aledom \disa, D_\my \aledom \mpsiu  )_{\aledom \Omega_f} + ( D_\my \aledom \dis, D_\my \aledom \mpsiw  )_{\aledom \Omega_f} \\
&\hspace{12pt} -  ( \aledom \disa, \aledom \mpsiw  )_{\aledom \Omega_f} +  ( \mathrm{div}_\my  ( {\maledomJ} \mFXa^{-1} \aledom \mv ), \aledom \mpsip  )_{\aledom \Omega_f} + (\aledom p, \aledom \mpsip)_{\aledom \Omega_s}  +  ( D_\my \aledom \disa, D_\my \aledom \mpsiw  )_{\aledom \Omega_s}   = 0,
\end{align*}
for all $ ( \aledom \mpsiv, \aledom \mpsip, \aledom \mpsiu, \aledom \mpsiw ) \in \aledom{\mathbf V}_0 \times \aledom P \times \aledom{\mathbf W}_{f,s,0}  \times \aledom{\mathbf Z}_{f,s,0}$ and a.e. $\mt \in I$.
For the sake of simplicity and in order to be able to work with standard finite element function spaces for the test functions in the numerical realization, we consider a modified weak formulation by replacing $ \aledom {\mathbf W}_{f,s,0}$ with $ \aledom {\mathbf W}_0$ and $ \aledom {\mathbf Z}_{f,s,0}$ with $ \aledom {\mathbf Z}$. This change corresponds to the introduction of additional interface integrals $- ( (D_\my \aledom \disa)\vert_{\aledom \Omega_s} \malens, \aledom \mpsiw)_{\aledom \Gamma_i}$ and $ - ( (D_\my \aledom \disa)\vert_{\aledom \Omega_f} \malenf, \aledom \mpsiu)_{\aledom \Gamma_i}$, which can be seen by performing an integration by parts:\\
Find $ ( \aledom \mv, \aledom p, \aledom \dis, \aledom \disa  ) \in W_{2,2}  (I, \aledom{\mathbf V}  ) \times L^2  (I, \aledom P   ) \times W_{2,2} ( I,  \aledom{\mathbf W}  )\times L^2  (I, \aledom{\mathbf Z}  )$ such that $\aledom \mv  \mevalo= \aledom \mv_0$, $\aledom \dis  \mevalo= \aledom \dis_0$ and
\begin{align*}
\begin{split}
&\langle \aledom A_3  ( \aledom \mv, \aledom p, \aledom \dis, \aledom \disa  ),  ( \aledom \mpsiv, \aledom \mpsip, \aledom \mpsiu, \aledom \mpsiw )  \rangle :=
( {\maledomJ} \rho_f \partial_\mt \aledom \mv, \aledom \mpsiv )_{\aledom \Omega_f} \\
&\hspace{12pt}+  ({\maledomJ} \rho_f  (  (\mFXa^{-1}  ( \aledom \mv - \partial_\mt \aledom \dis  ) )\cdot \nabla_\my  ) \aledom \mv, \aledom \mpsiv  )_{\aledom \Omega_f}  +  ( {\maledomJ} \aledom{ \boldsymbol \sigma }_f \mFXa^{-\top}, D_\my \aledom \mpsiv  )_{\aledom \Omega_f} -  (\aledom{\mathbf g}_f, \aledom \mpsiv  )_{\aledom \Gamma_{fN}} \\
&\hspace{12pt}-  ({\maledomJ} \rho_f \aledom{\mathbf f}_f, \aledom \mpsiv  )_{\aledom \Omega_f} +  ( \rho_s \partial_\mt \aledom \mv, \aledom \mpsiv  )_{\aledom \Omega_s} +  ( {\maledomJ} \aledom{\boldsymbol \sigma}_s \mFXa^{-\top}, D_\my \aledom \mpsiv  )_{\aledom \Omega_s} -  ( \aledom{\mathbf g}_s, \aledom \mpsiv  )_{\aledom \Gamma_{sN}} \\
&\hspace{12pt}-  ({\maledomJ} \rho_s \aledom{\mathbf f}_s, \aledom \mpsiv  )_{\aledom \Omega_s} +  ( \rho_s ( \partial_\mt \aledom \dis - \aledom \mv  ), \aledom \mpsiu  )_{\aledom \Omega_s} +  ( D_\my \aledom \disa, D_\my \aledom \mpsiu  )_{\aledom \Omega_f} + ( D_\my \aledom \dis, D_\my \aledom \mpsiw  )_{\aledom \Omega_f} \\
&\hspace{12pt} -  ( \aledom \disa, \aledom \mpsiw  )_{\aledom \Omega_f} +  ( \mathrm{div}_\my  ( {\maledomJ} \mFXa^{-1} \aledom \mv ), \aledom \mpsip  )_{\aledom \Omega_f} + (\aledom p, \aledom \mpsip)_{\aledom \Omega_s} +  ( D_\my \aledom \disa, D_\my \aledom \mpsiw  )_{\aledom \Omega_s}  \\
&\hspace{12pt}  - ( (D_\my \aledom \disa)\vert_{\aledom \Omega_f} \malenf, \aledom \mpsiu)_{\aledom \Gamma_i} - ( (D_\my \aledom \disa)\vert_{\aledom \Omega_s} \malens, \aledom \mpsiw)_{\aledom \Gamma_i}  = 0,
\end{split}
\end{align*}
for all $ ( \aledom \mpsiv, \aledom \mpsip, \aledom \mpsiu, \aledom \mpsiw ) \in \aledom{\mathbf V}_0 \times \aledom P \times \aledom{\mathbf W}_0 \times \aledom{\mathbf Z}$ and a.e. $\mt \in I$. 
Working with standard finite element spaces on the whole domain to discretize the above weak formulation results in a coupling between the original equations and the added auxiliary equations since the FEM functions cannot be arbitrarily localized (in fact, for continuous FEM nodal basis functions with a support intersecting both $\aledom \Omega_f$ and $\aledom \Omega_s$, the restriction to either $\aledom \Omega_f$ or $\aledom \Omega_s$ is not contained in the FEM space).
Hence, we introduce weighting parameters $\alpha_p>0$, $\alpha_w >0$, $\alpha_z > 0$ and functions $\aledom \iota_\bullet \in H^1(\aledom \Omega_\bullet)$, which are $0$ at the interface $\aledom \Gamma_i$, $1$ on $\aledom \Omega_\bullet \setminus B_\delta(\aledom \Gamma_i)$ and in $(0,1]$ on $B_\delta(\aledom \Gamma_i) \cap \aledom \Omega_\bullet$, where $\delta > 0$ is sufficiently small and $\bullet$ stands for either subscript $f$ or $s$.

The weak formulation, which we consider in the numerical implementation is given by: \\
Find $ ( \aledom \mv, \aledom p, \aledom \dis, \aledom \disa  ) \in W_{2,2}  (I, \aledom{\mathbf V}  ) \times L^2  (I, \aledom P   ) \times W_{2,2} ( I,  \aledom{\mathbf W}  )\times L^2  (I, \aledom{\mathbf Z}  )$ such that $\aledom \mv  \mevalo= \aledom \mv_0$, $\aledom \dis  \mevalo= \aledom \dis_0$ and
\begin{align}
\begin{split}
&\langle \aledom A  ( \aledom \mv, \aledom p, \aledom \dis, \aledom \disa  ),  ( \aledom \mpsiv, \aledom \mpsip, \aledom \mpsiu, \aledom \mpsiw )  \rangle :=
( {\maledomJ} \rho_f \partial_\mt \aledom \mv, \aledom \mpsiv )_{\aledom \Omega_f} \\
&\hspace{2pt}+  ({\maledomJ} \rho_f  (  (\mFXa^{-1}  ( \aledom \mv - \partial_\mt \aledom \dis  ) )\cdot \nabla_\my  ) \aledom \mv, \aledom \mpsiv  )_{\aledom \Omega_f}  +  ( {\maledomJ} \aledom{ \boldsymbol \sigma }_f \mFXa^{-\top}, D_\my \aledom \mpsiv  )_{\aledom \Omega_f} -  (\aledom{\mathbf g}_f, \aledom \mpsiv  )_{\aledom \Gamma_{fN}} \\
&\hspace{2pt}-  ({\maledomJ} \rho_f \aledom{\mathbf f}_f, \aledom \mpsiv  )_{\aledom \Omega_f} +  ( \rho_s \partial_\mt \aledom \mv, \aledom \mpsiv  )_{\aledom \Omega_s} +  ( {\maledomJ} \aledom{\boldsymbol \sigma}_s \mFXa^{-\top}, D_\my \aledom \mpsiv  )_{\aledom \Omega_s} -  ( \aledom{\mathbf g}_s, \aledom \mpsiv  )_{\aledom \Gamma_{sN}} \\
&\hspace{2pt}-  ({\maledomJ} \rho_s \aledom{\mathbf f}_s, \aledom \mpsiv  )_{\aledom \Omega_s} +  ( \rho_s ( \partial_\mt \aledom \dis - \aledom \mv  ), \aledom \mpsiu  )_{\aledom \Omega_s} +  \alpha_w( D_\my \aledom \disa, D_\my (\aledom \iota_f \aledom \mpsiu ) )_{\aledom \Omega_f} + ( D_\my \aledom \dis, D_\my \aledom \mpsiw )_{\aledom \Omega_f} \\
&\hspace{2pt} -  ( \aledom \disa,  \aledom \mpsiw  )_{\aledom \Omega_f} +  ( \mathrm{div}_\my  ( {\maledomJ} \mFXa^{-1} \aledom \mv ), \aledom \mpsip  )_{\aledom \Omega_f} + \alpha_p(\aledom p, \aledom \mpsip)_{\aledom \Omega_s} +  \alpha_z( D_\my \aledom \disa, D_\my ( \aledom \iota_s \aledom \mpsiw ) )_{\aledom \Omega_s}   = 0,
\end{split}
\label{eq::solopE}
\end{align}
for all $ ( \aledom \mpsiv, \aledom \mpsip, \aledom \mpsiu, \aledom \mpsiw ) \in \aledom{\mathbf V}_0 \times \aledom P \times \aledom{\mathbf W}_0 \times \aledom{\mathbf Z}$ and a.e. $\mt \in I$.

The corresponding formulation on the space-time cylinder reads as follows:\\
Find $ ( \aledom \mv, \aledom p, \aledom \dis, \aledom \disa  ) \in W_{2,2}  (I, \aledom{ \mathbf V}  ) \times L^2  (I, \aledom P   ) \times W_{2,2} ( I,  \aledom{ \mathbf W}  )\times L^2  (I, \aledom{\mathbf Z}  )$ such that $\aledom \mv  \mevalo= \aledom \mv_0$, $\aledom \dis  \mevalo= \aledom \dis_0$ and
\begin{align*}
&\langle \aledom A_{Q^T}   ( \aledom \mv, \aledom p, \aledom \dis, \aledom \disa  ),  ( \aledom \mpsiv, \aledom \mpsip, \aledom \mpsiu, \aledom \mpsiw )  \rangle :=
\mdl {\maledomJ} \rho_f \partial_\mt \aledom \mv, \aledom \mpsiv \mdr_{\mhQf} \\
& \hspace{2pt} +  \mdl {\maledomJ} \rho_f  (  (\mFXa^{-1}  ( \aledom \mv - \partial_\mt \aledom \dis  ) )\cdot \nabla_\my  ) \aledom \mv, \aledom \mpsiv  \mdr_{\mhQf} +  \mdl {\maledomJ} \aledom{ \boldsymbol \sigma }_f \mFXa^{-\top}, D_\my \aledom \mpsiv  \mdr_{\mhQf} -  \mdl \aledom{\mathbf g}_f, \aledom \mpsiv  \mdr_{\mhGfN} \\
& \hspace{2pt} -  \mdl{\maledomJ} \rho_f \aledom{\mathbf f}_f, \aledom \mpsiv  \mdr_{\mhQf} +  \mdl \rho_s \partial_\mt \aledom \mv, \aledom \mpsiv  \mdr_{\mhQs} +  \mdl {\maledomJ} \aledom{\boldsymbol \sigma}_s \mFXa^{-\top}, D_\my \aledom \mpsiv \mdr_{\mhQs} -  \mdl \aledom{\mathbf g}_s, \aledom \mpsiv  \mdr_{\mhGsN} \\
& \hspace{2pt} -  \mdl{\maledomJ} \rho_s \aledom{\mathbf f}_s, \aledom \mpsiv  \mdr_{\mhQs} +  \mdl \rho_s ( \partial_\mt \aledom \dis - \aledom  \mv  ), \aledom \mpsiu  \mdr_{\mhQs} + \alpha_w  \mdl D_\my \aledom \disa, D_\my (\aledom \iota_f \aledom \mpsiu)  \mdr_{\mhQf}  \\
& \hspace{2pt}+  \mdl D_\my \aledom \dis, D_\my \aledom \mpsiw  \mdr_{\mhQf}-  \mdl \aledom \disa, \aledom \mpsiw  \mdr_{\mhQf} + \alpha_p  \mdl \aledom p, \aledom \mpsip \mdr_{\mhQs} +  \mdl \mathrm{div}_\my  ( {\maledomJ} \mFXa^{-1} \aledom \mv ), \aledom \mpsip  \mdr_{\mhQf}  \\
&\hspace{2pt} +  \alpha_z \mdl D_\my \aledom \disa, D_\my ( \aledom \iota_s \aledom \mpsiw)  \mdr_{\mhQs} 
= 0,
\end{align*}
for all $ ( \aledom \mpsiv, \aledom \mpsip, \aledom \mpsiu, \aledom \mpsiw ) \in W_{2,2}  (I, \aledom{\mathbf V}_0  ) \times L^2  (I, \aledom P   ) \times W_{2,2} ( I,  \aledom{\mathbf W}_0  )\times L^2  (I, \aledom{\mathbf Z}  )$.

\section{Model of Shape Optimization Problem for FSI}
\label{sec::fsisop}
In this section we model the shape optimization problem (Section \ref{subsec::op}) for the unsteady, nonlinear FSI system \eqref{eq::solopE} via the method of mappings. For this purpose, we transform the FSI equations to the nominal domain (Section \ref{subsec::transform}), choose a set of admissible transformations (Section \ref{subsec::adtrans}) and an objective function (Section \ref{subsec::obfunc}).

\subsection{Transformation of FSI Equations to Nominal Domain}
\label{subsec::transform}
The method of mappings is applied to shape optimization problems that are governed by the FSI equations formulated on the ALE reference domain $\aledom \Omega$. The actual physical domain is obtained from the reference domain via the homeomorphism $\aledom{\mX}(\cdot, {\sf{t}}):  \aledom \Omega \to \pdom \Omega({\sf t})$. The method of mappings applies an additional transformation $\mmtau: \tilde \Omega \to \aledom\Omega$, which is a bi-Lipschitz transformation from the nominal domain $\tilde \Omega$ to $\aledom\Omega$. Let the design part of the
boundary $\sdom \Gamma_d$ be a subset of $(\partial \sdom \Omega \cap \partial \sdom \Omega_f) \cup \sdom \Gamma_i$.

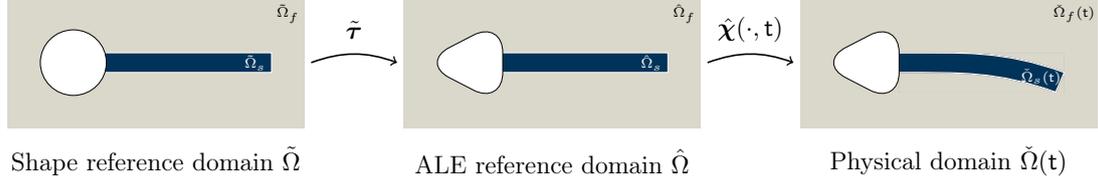
\begin{figure}[!h]
	\centering
	\resizebox{1.\linewidth}{!}{
		\begin{tikzpicture}[path A/.style 2 args={insert path={(13.5, 1.15) to #1 (13.5,0.85) #2}},
			path B/.style 2 args={insert path={(13.5, 1.15) to[out = 0, in = 160] #1 (16,0.85) #2}},
			path C/.style 2 args={insert path={((13.5, 0.85) to[out = 0, in = 160] #1 (15.87,0.55) #2}},
			path D/.style 2 args={insert path={(16,0.85) to #1 (15.87, 0.55) #2}}]
			\fill[tumg] (0.0,0.0) rectangle (4.5, 2.0);
			\fill[tumblues2, draw=white] (1.0,0.85) rectangle(4,1.15);
			\fill[white, draw = black] (1.0, 1.0) circle (0.5);
			\node[color = white] at (3.75, 1.0) {{\tiny{$\sdom\Omega_s$}}};
			\node at (4.25, 1.75) {{\tiny{$\sdom\Omega_f$}}};
			\node at (2.25, -0.5){Shape reference domain $\sdom\Omega$};
			\fill[tumg] (6.0,0.0) rectangle (10.5, 2.0);
			\fill[tumblues2, draw=white] (7.0,0.85) rectangle(10,1.15);		
			\draw [rounded corners=5mm,fill=white] (6.25,1.0)--(7.5,0.35)--(7.5,1.65)--cycle;
			\node[color = white] at (9.75, 1.0) {{\tiny{$\aledom\Omega_s$}}};
			\node at (10.25, 1.75) {{\tiny{$\aledom\Omega_f$}}};
			\node at (8.25,-0.5){ALE reference domain $\aledom\Omega$};
			\begin{pgfonlayer}{bg1}
				\fill[tumg] (12.0,0.0) rectangle (16.5, 2.0);
				\fill[tumblues2] (13.5,0.55) rectangle(16,1.15);
			\end{pgfonlayer}
			\draw[path A={coordinate[pos=0] (A1)
				coordinate[pos=1] (A2)}{}];
			\draw[path B={coordinate[pos=0] (B1) 
				coordinate[pos=1] (B2)}{}];
			\draw[path C={coordinate[pos=0] (C1)
				coordinate[pos=1] (C2)}{}];
			\draw[path D={coordinate[pos=0] (D1)
				coordinate[pos=1] (D2)}{}];
			\begin{pgfonlayer}{bg2}
				\clip[path C={}{-- (13.5,0.55) -- cycle}];
				\fill[tumg, draw =tumg] (13.5,1.15) -- (16,1.15) -- (16,0.55) -- (13.5,0.55) -- cycle;
			\end{pgfonlayer}
			\begin{pgfonlayer}{bg3}
				\clip[path B={}{-- (16,1.15) -- cycle}];
				\fill[tumg, draw = tumg] (13.5,1.15) -- (16,1.15) -- (16,0.55) -- (13.5,0.55) -- cycle;
			\end{pgfonlayer}
			\begin{pgfonlayer}{bg4}
				\clip[path D={}{-- (16,0.55) -- cycle}];
				\fill[tumg, draw = tumg] (13.5,1.15) -- (16,1.15) -- (16,0.55) -- (13.5,0.55) -- cycle;
			\end{pgfonlayer}
			\draw[color = tumg, line width = 0.4] (13.5, 1.15) to (16,1.15);
			\draw[color = tumg, line width = 0.4] (13.5, 0.55) to (16,0.55);
			\draw[color = tumg, line width = 0.4] (13.5, 1.15) to (13.5,0.55);
			\draw[color = tumg, line width = 0.4] (16, 1.15) to (16,0.55);
			\draw[color = white, line width = 0.5] (13.5, 1.15) to[out = 0, in = 160] (16,0.85);
			\draw[color = white, line width = 0.5] (13.5, 0.85) to[out = 0, in = 160] (15.87,0.55);
			\draw[color = white, line width = 0.5] (16, 0.85) to (15.87,0.55);
			\node[color = white] at (15.65, 0.785) {{\tiny{$\pdom\Omega_s({\sf t})$}}};
			\node at (16.15, 1.75) {{\tiny{$\pdom\Omega_f ({\sf t})$}}};
			\node at (14.25,-0.5){Physical domain $\pdom \Omega ({\sf t} )$};
			\draw [rounded corners=5mm,fill=white] (12.25,1.0)--(13.5,0.35)--(13.5,1.65)--cycle;
			\draw[->, line width = 0.7pt] (4.6,1.0) to[out=20, in=160] (5.9,1.0);
			\draw[->, line width = 0.7pt] (10.6,1.0) to[out=20, in=160] (11.9,1.0);
			\node at (5.25,1.5) {$\mmtau$};
			\node at (11.25,1.5){${\aledom \mX} (\cdot,{\sf t })$};
		\end{tikzpicture}
	}
	\caption{Schematic illustration of the method of mappings combined with an ALE transformation}
	\label{fsi::fig::trafo}
\end{figure}

The shape reference, ALE, and physical domains and the transformations between them are illustrated in Figure \ref{fsi::fig::trafo}. 

Since we want to optimize the shape of the domain and not the initial conditions or boundary conditions, we assume that the considered transformations in $\mTad$ do not change these conditions. For the sake of convenience, and in correspondence with our numerical setting, we choose $\aledom{\mathbf g}_f = 0$ and $\aledom{\mathbf g}_s = 0$. Additionally, we choose the transformation such that it is equal to the identity on the support of the initial conditions (which is not a restriction in our case since the initial conditions are chosen to be $0$).\\
For fixed $\mmtau \in \mTad$, we introduce the  following spaces on the shape reference domain
\begin{align*}
&\sdom{\mathbf V}_\mmtau = \lbrace \aledom \vel \circ \mmtau ~:~ \aledom \vel \in \aledom{\mathbf V}(\mmtau(\sdom \Omega)) \rbrace, \\
&\sdom{\mathbf V}_{0, \mmtau} = \lbrace \aledom \vel \circ \mmtau ~:~ \aledom \vel \in \aledom{\mathbf V}_0(\mmtau(\sdom \Omega)) \rbrace, \\
&\sdom{\mathbf W}_\mmtau = \lbrace \aledom \dis \circ \mmtau ~:~ \aledom \dis \in \aledom{\mathbf W}(\mmtau(\sdom \Omega)) \rbrace, \\
&\sdom{\mathbf W}_{0, \mmtau} = \lbrace \aledom \dis \circ \mmtau ~:~ \aledom \dis \in \aledom{\mathbf W}_0(\mmtau(\sdom \Omega)) \rbrace, \\
&\sdom P_\mmtau = \lbrace \aledom p \circ \mmtau ~:~ \aledom p \in \aledom P(\mmtau(\sdom \Omega)) \rbrace, \\
&\sdom{\mathbf Z}_\mmtau = \lbrace \aledom \disa \circ \mmtau ~:~ \aledom \disa \in \aledom{\mathbf Z}(\mmtau(\sdom \Omega)) \rbrace.
\end{align*}
\begin{remark}
    For a function $\aledom \qfct: ~ \aledom \Omega \times I \to \mathbb R^n$, $n \geq 1$, $\aledom \qfct \circ \mmtau$ denotes the function that maps $(\mz, \mt) \mapsto \aledom \qfct(\mmtau(\mz), \mt)$. 
\end{remark}
The additional transformation with $\mmtau$ yields the following weak formulation on the shape reference domain $\sdom \Omega$. \\
For fixed $\mmtau \in \mTad$, find 
$ ( \sdom  \mv, \sdom  p, \sdom \dis, \sdom \disa) \in W_{2,2}(I, \sdom{\mathbf V}_\mmtau) \times L^2(I, \sdom P_\mmtau) \times W_{2,2}(I, \sdom{\mathbf W }_\mmtau) \times L^2(I, \sdom{\mathbf Z}_\mmtau ) $ such that $\sdom \mv\mevalo = 0$, $ \sdom \dis\mevalo = 0$ and
\begin{align}
\begin{split}
& \langle \sdom A_{ \Omega}  ( (\sdom \mv, \sdom p, \sdom \dis, \sdom \disa ), \mmtau ),  ( \sdom \mpsiv, \sdom \mpsip, \sdom \mpsiu, \sdom \mpsiw )  \rangle \\&:=
(\mathrm{det}(D_\mz \mmtau) {\msdomJ} \rho_f \partial_\mt \sdom \mv, \sdom \mpsiv )_{\sdom \Omega_f} \\ & \hspace{12pt}+  (\mathrm{det}(D_\mz \mmtau) {\msdomJ} \rho_f  (  ((D_\mz \mmtau)^{-1}{\msdomJXt}^{-1}  ( \sdom \mv - \partial_\mt \sdom \dis  ) )\cdot \nabla_\mz  ) \sdom \mv, \sdom \mpsiv  )_{\sdom \Omega_f} \\
& \hspace{12pt}+  ( \mathrm{det}(D_\mz \mmtau){\msdomJ} \sdom{ \boldsymbol \sigma }_f {\msdomJXt}^{-\top}(D_\mz \mmtau)^{-\top}, D_\mz \sdom \mpsiv  )_{\sdom \Omega_f}  -  (\mathrm{det}(D_\mz \mmtau){\msdomJ} \rho_f \sdom{\mathbf f}_f, \sdom \mpsiv  )_{\sdom \Omega_f} \\ 
&\hspace{12pt}+  (\mathrm{det}(D_\mz \mmtau) \rho_s \partial_\mt \sdom \mv, \sdom \mpsiv  )_{\sdom \Omega_s}  +  (\mathrm{det}(D_\mz \mmtau) {\msdomJ} \sdom{\boldsymbol \sigma}_s {\msdomJXt}^{-\top} (D_\mz \mmtau)^{-\top}, D_\mz \sdom \mpsiv  )_{\sdom \Omega_s}\\
& \hspace{12pt} -  (\mathrm{det}(D_\mz \mmtau){\msdomJ} \rho_s \sdom{\mathbf f}_s, \sdom \mpsiv  )_{\sdom \Omega_s}  +  ( \mathrm{det}(D_\mz \mmtau) \rho_s ( \partial_\mt \sdom \dis - \sdom \mv  ), \sdom \mpsiu  )_{\sdom \Omega_s} \\
&\hspace{12pt}+ \alpha_w  ( \mathrm{det}(D_\mz \mmtau) D_\mz \sdom \disa (D_\mz \mmtau)^{-1} (D_\mz \mmtau)^{-\top}, D_\mz (\sdom \iota_f \sdom \mpsiu)  )_{\sdom \Omega_f} \\
&\hspace{12pt} +  ( \mathrm{det}(D_\mz \mmtau) D_\mz \sdom \dis (D_\mz \mmtau)^{-1} (D_\mz \mmtau)^{-\top} , D_\mz \sdom \mpsiw  )_{\sdom \Omega_f} \\
&\hspace{12pt} -  (\mathrm{det}(D_\mz \mmtau) \sdom \disa, \sdom \mpsiw  )_{\sdom \Omega_f} + \alpha_p  (\mathrm{det}(D_\mz \mmtau) \sdom p,  \sdom \mpsip )_{\sdom \Omega_s} \\ 
&\hspace{12pt}+  ( \mathrm{det}(D_\mz \mmtau) \mathrm{tr}( D_\mz ( {\msdomJ} {\msdomJXt}^{-1} \sdom \mv ) (D_\mz \mmtau)^{-1}), \sdom \mpsip  )_{\sdom \Omega_f} \\
&\hspace{12pt} +  \alpha_z( \mathrm{det}(D_\mz \mmtau) D_\mz \sdom \disa (D_\mz \mmtau)^{-1} (D_\mz \mmtau)^{-\top}, D_\mz (\sdom \iota_s \sdom \mpsiw)  )_{\sdom \Omega_s}  = 0,
\end{split}
\label{eq::transformedfsi}
\end{align}
for all $ ( \sdom \mpsiv, \sdom \mpsip, \sdom \mpsiu, \sdom \mpsiw ) \in W_{2,2}  (I, \sdom{\mathbf V}_0  ) \times L^2  (I, \sdom P   ) \times W_{2,2} ( I,  \sdom{\mathbf W}_0  )\times L^2  (I, \sdom{\mathbf Z}  )$ and any $\mt \in I$. Moreover, ${\msdomJ} = \mathrm{det} ({\msdomJXt})$ and ${\msdomJXt} = D_\mz \sdom \mX  (D_\mz \mmtau)^{-1}$,
\[
\sdom{ \boldsymbol \sigma}_{f} := \rho_f \nu_f  ( D_\mz \sdom{\mv}_f (D_{\mz} \mmtau)^{-1}{\msdomJXt}^{-1} + {\msdomJXt}^{-\top} (D_{\mz} \mmtau)^{-\top} D_\mz \sdom{ \mv}_f^{\top}  ) - \sdom{p}_f \mathbf \mId 
\]
denotes the transformed fluid stress tensor and $\sdom \sigma_s$ the corresponding transformed solid stress tensor. For Saint Venant-Kirchhoff type material, it is given by
\[
\sdom{\boldsymbol \sigma}_s = \msdomJ^{-1}{\msdomJXt} ( \lambda_s \mathrm{tr} ({\msdomE}  )\mathbf I + 2 \mu_s {\msdomE}) {\msdomJXt^\top}
\]
with ${\msdomE}  := \frac{1}{2} ({\msdomJXt}^{\top} {\msdomJXt}  - \mId  )$. The corresponding operator is denoted by 
\begin{align} \sdom A_{\Omega}( \sdom{\mathbf y}, \mmtau ) = 0,\label{equation::fsiop}\end{align} where $\sdom{\mathbf y} = (\sdom \mv, \sdom p, \sdom \dis, \sdom \disa)$. Moreover, let $$\sdom{\mathbf Y}_\mmtau := W_{2,2}(I, \sdom{\mathbf V}_\mmtau) \times L^2(I, \sdom P_\mmtau) \times W_{2,2}(I, \sdom{\mathbf W }_\mmtau) \times L^2(I, \sdom{\mathbf Z}_\mmtau ).$$ 

\subsection{Choice of Objective Function}
\label{subsec::obfunc}
As objective function, we choose the mean fluid drag which is given by 
\[
- \frac1T \int_0^T \int_{ \check \Gamma_o(\mt)} \pdom {\boldsymbol \psi}^{\top} \boldsymbol{\sigma}_{f,\mx} (\pdom \mv_f, \pdom \p_f) \pdom{\mathbf n}_f dS(\mx) d\mt 
\]
where $\check \Gamma_o(\mt) \subset \partial \check\Omega_f (\mt)$ for all $\mt \in [0,T]$
denotes an obstacle such that $\check \Gamma_o(\mt)$ and $\partial \check\Omega_f (\mt) \setminus \check \Gamma_o(\mt)$ have positive distance for all $\mt \in [0,T]$.
Furthermore, 
$\pdom {\boldsymbol \psi} = (1, 0)^{\top}$ and $ \pdom{\mathbf n}_f$ denotes the outwards pointing normal vector, e.g., 
\cite{BaLiUl, JoMa01}.
This can be reformulated as a volume integral given by
\[
-\frac1T \int_0^T ( ( \rho_f(\partial_\mt \pdom \mv_f + (\pdom \mv_f \cdot \nabla_\mx) \pdom \mv_f ), \pdom {\boldsymbol \Psi})_{\pdom \Omega_f(\mt)} - (\pdom \p_f, \mathrm{div}(\pdom {\boldsymbol \Psi}))_{\pdom \Omega_f(\mt)} + (2 \nu_f \epsilon_\mx(\pdom \mv_f), \epsilon_\mx(\pdom {\boldsymbol \Psi}))_{\pdom \Omega_f (\mt)} ) d\mt,
\]
where $\epsilon_\mx(\cdot) = \frac12 (D_\mx \cdot + (D_\mx \cdot)^{\top})$, and $\pdom {\boldsymbol \Psi}$ is an arbitrary continuously differentiable function such that $\pdom {\boldsymbol \Psi} \vert_{\check \Gamma_o} = \pdom{\boldsymbol \psi}$ and $\pdom  {\boldsymbol \Psi} \vert_{\partial {\check \Omega_f}\setminus \check \Gamma_o} = 0$, cf. 
\cite{BraRi, HoJo02, JoMa01}
. 
\begin{remark}
	Analogously to the observation in \cite{HPS15}, where the evaluation of the shape gradient via surface integrals is compared to the evaluation via volume integrals, using the volume integral formulation for the drag is expected to provide more accurate numerical results. Our numerical tests confirmed this.
\end{remark}
The corresponding transformed formulation on the ALE domain $\aledom \Omega_f$ reads as
\begin{align*}
\widealedom{ F}_D (\aledom {\boldsymbol y})= & - \frac1T \int_0^T ( ( {\maledomJ} \rho_f(\partial_\mt \aledom \mv_f+(({\maledomJX}^{-1}(\aledom \mv_f - \partial_\mt \aledom \mX))\cdot \nabla_\my)\aledom \mv_f), \aledom {\boldsymbol \Psi}  )_{\aledom \Omega_f} \\ &- ({\maledomJ} \aledom p_f, \mathrm{tr}(D_\my \aledom {\boldsymbol \Psi} {\maledomJX}^{-1} )  )_{\aledom \Omega_f} + (2\nu_f {\maledomJ} \epsilon_\my (\aledom \mv_f), \epsilon_\my (\aledom {\boldsymbol \Psi}))_{\aledom \Omega_f} ) d\mt
\end{align*}
with $ {\maledomJ} = \mathrm{det}{\maledomJX}$ and $\epsilon_\my (\cdot) = \frac12 ((\mDy \cdot){\maledomJX}^{-1} + {\maledomJX}^{-\top} (\mDy \cdot )^\top )$, and the transformation on the shape reference domain $\sdom \Omega_f$ yields
\begin{align*}
\sdom{ F}_D (\sdom{\mathbf y}, \mmtau) = & - \frac1T \int_0^T (( {\msdomJ} \mathrm{det} (D_\mz \mmtau) \rho_f (\partial_\mt \sdom \mv_f + (( D_\mz \mmtau^{-1} {\msdomJXt}^{-1}(\sdom \mv_f - \partial_\mt \sdom \mX))\cdot \nabla_\mz)\sdom \mv_f), \sdom {\boldsymbol \Psi}    )_{\sdom \Omega_f} \\ &
- ({\msdomJ} \mathrm{det} (D_\mz \mmtau)  \sdom p, \mathrm{tr}(D_\mz \sdom {\boldsymbol \Psi} D_\mz \mmtau^{-1} {\msdomJXt}^{-1}  ))_{\sdom \Omega_f} \\ & + (2 \nu_f {\msdomJ} \mathrm{det} (D_\mz \mmtau) \epsilon_\mz (\sdom \mv_f),   \epsilon_\mz (\sdom {\boldsymbol \Psi}) )_{\sdom \Omega_f} )d\mt
\end{align*}
with $\epsilon_\mz (\cdot) = 
\frac12( (\mDz \cdot) \mDz \mmtau^{-1} \msdomJXt^{-1} + \msdomJXt^{-\top}
\mDz \mmtau^{-\top} (\mDz \cdot)^\top )$.

\subsection{Choice of Admissible Shape Transformations}
\label{subsec::adtrans}
Our choice of admissible transformations is based on the following considerations:
\begin{itemize}
	\item As already mentioned in a previous section, we choose the transformations such that they do not change initial conditions, boundary conditions or source terms, i.e., the support of the deformation $\mmtau - \mathrm{id}_\mz$ is disjoint from the support of the initial conditions, boundary conditions and source terms. In case that the design part $\sdom \Gamma_d$ is a subset of the Dirichlet boundary, the boundary conditions are considered to be homogeneous on $\sdom \Gamma_d$.   
	\item Since standard existence theory for PDEs requires Lipschitz regularity of the domain, it is straightforward to require the domains to be Lipschitzian during the optimization process. This can be ensured by choosing $\sdom \Omega$ as a Lipschitz domain and transformations $\mmtau \in W^{1, \infty}(\sdom \Omega)^d$ close to the identity \cite[Lem. 2]{BeFe}. The locality around the identity can be relaxed by including a determinant constraint $\mathrm{det}(D \mmtau) \geq \eta_{\mathrm{ext}}$, $\eta_{\mathrm{ext}} \in (0,1)$, on $\Omega$ or an hold all domain $U \supset \Omega$, see \cite{HSU20, etling2020first}. In this work, we choose this perspective and add a penalization $ \gamma_p \int_{\tilde \Omega} ( \mathrm{det}(D_\mz \mmtau(\mz)) - \eta_{\mathrm{ext}})^{-1} d \mz$ to the objective function, with $\gamma_p = 10^{-3}$ and $\eta_{\mathrm{ext}} = 0.2$, and set the value of the integral to $+ \infty$ (or in the numerical discretization to $10^{16}$) if the constraint $\mathrm{det}(D_\mz \mmtau) - \eta_{\mathrm{ext}} \geq 0$ is not fulfilled a.e. 
	\item To be able to work with the same function spaces on the shape reference domain independently of the control $\mmtau$, it is desirable that the spaces for the transformed functions are isomorphic to the spaces on the transformed domain.
    When working in a space $\sdom{\mathbf Y} = \sdom{\mathbf Y}(\sdom \Omega)$ on the shape reference domain, the corresponding space on the ALE reference domain is $$\aledom{\mathbf Y}_\mmtau (\aledom \Omega) = \lbrace \aledom{\mathbf y}~:~ \aledom{\mathbf y} \circ \mmtau \in \sdom{\mathbf Y} \rbrace.$$
    Hence, the regularity of the space on the actual (in our case ALE reference) domain $\aledom \Omega$ depends on the regularity of $\mmtau$ and also the regularity of $\aledom \Omega$ depends on the regularity of $\mmtau$ and $\sdom \Omega$.
    One way to ensure well-posedness of the optimization problem is the choice of a suitable notion of solutions on $\aledom \Omega$, a function space setting on $\aledom \Omega$, and $\mmtau$ such that solving the PDE on the actual domain is equivalent to solving the transformed PDE on the shape reference domain and $\aledom {\mathbf Y} = \aledom{\mathbf Y}_\mmtau (\aledom \Omega)$ is independent of $\mmtau$ \cite{HUU19}.
    Thus, the regularity requirement on $\mmtau$ depends on the regularity of the state of the PDEs and $\mTad \subset  \sdom {\mathbf D}_{ \Omega} \subset W^{1,\infty}(\sdom \Omega)^d$ for a function space $\sdom {\mathbf D}_{ \Omega}$ with sufficiently high regularity. 
	\item Transformations that only change the interior of the domain but not the boundaries do not change the shape of the domain. To ensure a one-to-one correspondence, shape optimization problems are often considered as optimization problems on manifolds, see, e.g., 
    \cite{MiMu, Ring}, or on appropriate subsets of linear subspaces, see, e.g., \cite{BaLiUl}.
	In order to be in the latter setting, we consider a scalar valued quantitiy $\sdom d \in  \sdom D_{ \Gamma_d}$ on the design boundary $\sdom \Gamma_d$ and identify it with a shape via a transformation of the form $\mathrm{id}_\mz + \mathrm B(\sdom d)$. 
\end{itemize}

\noindent
We work with the following operator $B: L^2(\sdom \Gamma_d) \to \sdom{\mathbf D}_\Omega, ~ \sdom d \mapsto \sdom{\mathbf u}$: We choose $ \sdom D_{\Gamma_d} = L^2(\sdom \Gamma_d)$ and 
	first solve the Laplace-Beltrami equation 
	\begin{equation}\label{LapBel}
	 - \mathrm{div}_{\sdom \Gamma_d} ( \tilde \beta D_{\sdom \Gamma_d} \sdom{\mathbf b} ) + \sdom{\mathbf b} = \sdom d \sdom{ \mathbf{n}} ~ \text{in }\sdom \Gamma_d, \quad  \sdom{\mathbf b} = 0 ~ \text{on } \partial \sdom \Gamma_d, 
	\end{equation}
	where $\tilde \beta$ is an $H^1$-regularized approximation to $\beta>0$ in $L^2$ that is $1$ on $\partial \sdom \Gamma_d$. More precisely, $\tilde \beta$ solves
	\begin{equation}\label{tildebetaeq}
		- \beta \Delta_{\sdom \Gamma_d} \tilde \beta + \tilde \beta = \beta  ~ \text{in } \sdom \Gamma_d, \quad \tilde \beta = 1 ~ \text{on } \partial \sdom \Gamma_d,
	\end{equation}
	where $\beta >0$ is given. In contrast to \cite{HSU20}, $\tilde \beta$ is chosen via a PDE on $\tilde \Gamma_d$ since we consider design boundaries which have a boundary themselves. In order to avoid mesh degeneration, we aim for small gradients of $\sdom{\mathbf b}$ in the vicinity of $\partial \sdom \Gamma_d$. 
	To obtain the deformation field, we solve the equation
	\begin{align}
    \begin{split}
		&- \mathrm{div}_{\sf z}(D_{\sf z} \sdom{\mathbf u}_\tau + D_{\sf z} \sdom{\mathbf u}_\tau^\top ) = 0 ~ \text{in } \Omega, \quad \sdom{\mathbf u}_\tau = 0 ~\text{on } (\partial \sdom \Omega \cup \sdom \Gamma_i) \setminus \sdom \Gamma_d, \\ & (D_{\sf z} \sdom{\mathbf u}_\tau + D_{\sf z} \sdom{\mathbf u}_\tau^\top) \sdom{\mathbf n} = \sdom{\mathbf b} ~\text{on } \sdom \Gamma_d,
    \end{split}
    \label{extension::equation}
	\end{align}
	which is similarly also used in the traction method \cite{Az70} or in Steklov-Poincar\'e type methods \cite{schulz2016efficient}. Note that we work with (too) smooth transformations and can therefore not obtain additional kinks. Driving $\beta$ to $0$ allows for the approximation of kinks, see \cite{HSU20}.

\begin{remark}
The above considerations allow for a variety of possible choices for the operator $\mathrm B$. In \cite{HDiss} we work with a smooth deformation direction field that is appropriately scaled by a scalar valued function that depends on $\sdom d$. Here, we consider a strategy similar to the one introduced in \cite{HSU20} which is not tailored towards the geometrical setting but allows for generalization to different geometries and is also applicable in 3D.
\end{remark}

We consider sets of admissible transformations
\[
\mTad \subset \lbrace \mmtau = \mathrm{id}_\mz + \mmut,\, \mmut \in \mUt \rbrace, 
\]
where $\mUt$ is chosen such that 
\[
\mUt \subset \lbrace \mmut ~:~ \mmut= \mathrm B(\sdom d),\, \sdom d \in  \sdom D_{\Gamma_d} \rbrace.
\]
Here, $\mUt$ is chosen as a closed subset. In order to ensure the 
bi-Lipschitzian property for all transformations corresponding to an admissible control, additional conditions have
to be satisfied, e.g., $\mlnorm \sdom d \mrnorm_{ \sdom D_{\Gamma_d}} \leq c$ 
for a sufficiently small constant $c>0$, or a lower bound on the determinant of the gradient of the deformation.  \\[1ex]
Furthermore, it is often relevant for practical applications to have additional geometric constraints, e.g., that the volume and barycenter of the obstacle shall not change. This motivates the constraints
\begin{align*}
\sdom g_\Omega(\mmut) = \begin{pmatrix} 
&\int_{\sdom \Omega} \mathrm{det} (\mathbf I + D_\mz \mmut )  d \mz - V \\
& \frac1{V_o} (C_1 - \int_{{\sdom \Omega}} (\mz + \mmut)_1 \mathrm{det} (\mathbf I + D_\mz \mmut ) d \mz) - B_1 \\
& \vdots \\
& \frac1{V_o} (C_d - \int_{{\sdom \Omega}} (\mz + \mmut)_d \mathrm{det} (\mathbf I + D_\mz \mmut ) d \mz) - B_d
\end{pmatrix}= 0,
\end{align*}
where $\sdom B_o$ denotes the obstacle, $V_o = \int_{\sdom B_o} \mathrm{det} (\mathbf I + D_\mz \mmut ) d \mz$ is the volume of the obstacle, $B = (B_1, \ldots, B_d)$ denotes the barycenter of the obstacle, $C = (C_1, \ldots, C_d)$ is the barycenter of $\sdom B_o \cup \sdom \Omega$, and $V$ denotes the volume of $\sdom \Omega$. Moreover, we assumed that $\mmtau(\sdom \Omega) \subset (\sdom B_o \cup \sdom \Omega)$, the transformed obstacle is defined by $(\sdom B_o \cup \sdom \Omega) \setminus \mmtau(\sdom \Omega)$, and there exists $\mmtau_{ext}: \sdom B_o \cup \sdom \Omega \to \sdom B_o \cup \sdom \Omega$ that is bi-Lipschitz and such that $\mmtau_{ext}\vert_{\sdom \Omega} = \mmtau$. In the numerical implementation, we use that the volume conservation constraint implies $V_o = \int_{\sdom B_o} 1 d \mz$ and thus we use this simpler formula.
\label{sec::shapetrans}

\subsection{Shape Optimization Problem}
\label{subsec::op}
The shape optimization problem is given by 
\begin{align}
\begin{split}
\min_{\sdom d \in \sdom D_{\Gamma_d} } &\tilde j_{\Omega} (\mmut) + \mathcal R(\sdom d) \\
\mathrm{s.t. } \quad & \tilde g_\Omega (\mmut) = 0,\quad \mmut = \mathrm B(\tilde d),
\end{split}
\label{cont::opt::1}
\end{align}
where $\sdom j_{\Omega} (\mmut) = \sdom F_D (\sdom{\mathbf y },\mathrm{id}_\mz +  \mmut)$, and $\sdom {\mathbf y}$ is the solution to the PDE $\sdom A_{\Omega}( \sdom {\mathbf y}, \mathrm{id}_\mz + \mmut ) = 0$, see \eqref{equation::fsiop}. Furthermore, $\mathrm B$ and $\sdom g_\Omega$ are defined in Section \ref{sec::shapetrans} and $\sdom F_D$ is defined in Section \ref{subsec::obfunc}. Moreover, we choose $ \mathcal R (\sdom d) = \frac{\alpha}{2} \| \tilde d \|_{\sdom D_{\Gamma_d}}^2 + \mathcal R_p(\sdom d)$ with $\alpha > 0$ and $\mathcal R_p(\sdom d) = \gamma_p \int_{\tilde \Omega} ( \mathrm{det}(D_\mz \mmtau(\mz)) - \eta_{\mathrm{ext}})^{-1} d \mz$ being an additional penalization of a determinant constraint (see Section \ref{subsec::adtrans}). 

\section{Discretization}
\label{sec::fsidis}
In this section, we discretize the FSI system \eqref{eq::transformedfsi} in time (Section \ref{subsec::disct}) and space (Section \ref{subsec::disch}). To obtain a discrete formulation (Section \ref{subsec::dop}) of the optimization problem \eqref{cont::opt::1}, the objective function (Section \ref{subsec::discof}) and the shape transformations (Section \ref{subsec::discst}) have to be discretized. 

\subsection{Temporal Discretization}
\label{subsec::disct}
Our discretization in time uses a One-Step-$\theta$ scheme, cf. \cite{Wi11}. For this, we divide the terms that appear in the weak formulation into different categories. The first group $\sdom A_T (\tilde{\mathbf y}, \tau )  (\tilde{\boldsymbol \psi} )$ collects all terms which include time derivatives; further below, it then will be discretized in time by finite differences:
\begin{align*}
\sdom A_T (\sdom{\mathbf y}, \mmtau )  (\tilde{\boldsymbol \psi} ) :=\, & (\mtJtt \mtJht \rho_f  ( \partial_\mt \sdom \vel -   ( (\mtFti \mtFht^{-1} \partial_\mt \mtuht  ) \cdot {\nabla}_\mz  ) \sdom \vel) , \sdom \mpsiv  )_{\sdom{\Omega}_f} \\ & +  (\mtJtt {\rho}_s \partial_\mt \sdom \vel, \sdom \mpsiv  )_{\sdom{\Omega}_s} +  (\mtJtt  \rho_s \partial_\mt \mtuht , \sdom \mpsiu  )_{\sdom\Omega_s}.
\end{align*}
The group $\sdom A_I (\tilde{\mathbf y}, \tau )  (\tilde{\boldsymbol \psi} )$ gathers all implicit terms, i.e., all terms that should be fulfilled exactly by the new iterate such as the incompressibility condition for the fluid:
\begin{align*}
\sdom A_I (\sdom{\mathbf y}, \mmtau )  (\sdom{\boldsymbol \psi} ) :=\, &
(\mtJtt \mathrm{tr} (D_\mz  (\mtJht \mtFht^{-1}\sdom \vel ) \mtFti  )  , \sdom \mpsip )_{\sdom{\Omega}_f}  \\
& +  {\alpha}_p  ( \mtJtt \sdom p, \sdom \mpsip )_{\sdom \Omega_s}
-  ( \mtJtt \sdom \disa, \sdom \mpsiw )_{\sdom\Omega_f} \\ & +  (  \mtJtt D_\mz \mtuht \mtFti \mtFi ,  D_\mz \sdom \mpsiw )_{\sdom\Omega_f} \\
& +  \alpha_z( \mathrm{det}(D_\mz \mmtau) D_\mz \sdom \disa (D_\mz \mmtau)^{-1} (D_\mz \mmtau)^{-\top}, D_\mz (\sdom \iota_s \sdom \mpsiw)  )_{\sdom \Omega_s}.
\end{align*}
Another group $\sdom A_P (\sdom{\mathbf y}, \mmtau )  (\sdom{\boldsymbol \psi} )$, which is also treated implicitly, collects the pressure terms:
\begin{align*}
\sdom A_P (\sdom{\mathbf y}, \mmtau )  (\sdom{\boldsymbol \psi} ) :=  (\mtJtt \mtJht \tilde{\boldsymbol{\sigma}}_{f,p} \mtFht^{-\top} \mtFi, D_\mz \sdom \mpsiv   )_{\sdom{\Omega}_f},
\end{align*}
where $\tilde{\boldsymbol{\sigma}}_{f,p} = -\tilde p \mId$.
This can be motivated by the fact that the pressure serves as Lagrange multiplier for the incompressibility condition.
The remaining terms are collected in the fourth group $\sdom A_E (\tilde{\mathbf y}, \mmtau )  (\tilde{\boldsymbol \psi} )$:
\begin{align*}
\sdom A_E (\tilde{\mathbf y}, \mmtau )  (\tilde{\boldsymbol \psi} ) :=\,
&  (\mtJtt \mtJht \rho_f  ( (\mtFti \mtFht^{-1}  \sdom \vel  ) \cdot {\nabla}_\mz  ) \tilde{\vel}, \sdom \mpsiv  )_{\sdom{\Omega}_f} \\ &+ \alpha_w  ( \mtJtt D_\mz \sdom{\disa} \mtFti \mtFi,  D_\mz (\sdom \iota_f \sdom \mpsiu) )_{\sdom \Omega_f}   \\
&+  (\mtJtt \mtJht \tilde{\boldsymbol{\sigma}}_{f,v} \mtFht^{-\top} \mtFi, D_\mz \sdom \mpsiv   )_{\sdom{\Omega}_f} \\& +  (\mtJtt \mtJht \tilde{\boldsymbol{\sigma}}_s \mtFht^{-\top} \mtFi, D_\mz \sdom \mpsiv   )_{\sdom{\Omega}_s} \\ &-  ( \mtJtt \mtJht  \rho_f \tilde{\mathbf f}_f, \sdom \mpsiv  )_{\sdom {\Omega}_f}   -  ( \mtJtt \mtJht  \rho_s \tilde{\mathbf f}_s, \sdom \mpsiv )_{\sdom{\Omega}_s} \\&-  (\mtJtt  \rho_s \sdom \vel, \sdom \mpsiu  )_{\sdom\Omega_s},
\end{align*}
where $\sdom{\boldsymbol \sigma}_{f,v} = \sdom{\boldsymbol \sigma}_{f} - \sdom{\boldsymbol \sigma}_{f,p}$.
The time-stepping scheme can thus be summarized as follows.
Let a transformation $\mmtau$ be given, $N \in \mathbb N$, $0=t_0 < t_1 < \ldots < t_N = T$ be a discretization of $\overline I=[0,T]$ and $\theta \in  [0,1 ]$. Let, for $j \in  \lbrace 1,2,\ldots,N \rbrace$, $\tilde{\mathbf y}^{j-1}$ be the solution at the time $t_{j-1}$ and the time step size be constant, i.e., $k:= k_j = t_j - t_{j-1}$ for all $n \in \lbrace 1, \dots, N\rbrace$. Then, the solution at $t_j$ is computed by:\\
Find $\tilde{\mathbf y}^j \in \sdom{\mathbf V}_\mmtau \times \sdom P_\mmtau \times \sdom{\mathbf W }_\mmtau \times \sdom{\mathbf Z}_\mmtau$ such that 
\begin{align*}
&\sdom A_T^{j,k} (\sdom{\mathbf y}^j, \mmtau )  (\tilde{\boldsymbol \psi} ) + \theta \sdom A_E (\sdom{\mathbf y}^j, \mmtau )  (\sdom{\boldsymbol \psi} ) + \sdom A_P  (\sdom{\mathbf y}^j, \mmtau ) (\sdom{\boldsymbol\psi} ) + \sdom A_I (\sdom{\mathbf y}^j, \mmtau ) (\sdom{\boldsymbol\psi} ) \\ & 
= -  (1 -\theta ) \sdom A_E (\sdom{\mathbf y}^{j-1}, \mmtau ) (\sdom{\boldsymbol \psi} ),
\end{align*}
for all test functions $\sdom{\boldsymbol \psi} \in \sdom{\mathbf V}_0 \times \sdom P \times \sdom{\mathbf W}_0 \times \sdom{\mathbf Z}$.
Here, $ \sdom A_T^{j,k} (\sdom{\mathbf y},\mmtau ) (\sdom{\boldsymbol\psi} )$ is defined as the approximation of $\sdom A_T (\sdom{\mathbf y}, \mmtau)(\sdom{\boldsymbol\psi} )$ given by
\begin{align*}
&\sdom A_T^{j,k} (\sdom{\mathbf y},\mmtau ) (\sdom{\boldsymbol\psi} )\\&:= 
\frac{1}{k}  (\mtJtt \mtJht^{j,\theta} \rho_f  ( (\sdom{\vel} - \sdom {\vel}^{j-1}) -   ( (\mtFti \mtFht^{-1} ( \mtuht - \mtuht^{j-1}  ) \cdot {\nabla}_\mz  ) \sdom{\vel} ) , \sdom \mpsiv  )_{\sdom{\Omega}_f} \\ & \hspace{16pt}
+ \frac{1}{k}  (\mtJtt {\rho}_s (\sdom {\vel} - \sdom {\vel}^{j-1}), \sdom\mpsiv )_{\sdom{\Omega}_s} + \frac1k  (\mtJtt  \rho_s (\mtuht - \mtuht^{j-1}) , \sdom \mpsiu  )_{\sdom\Omega_s},
\end{align*}
where $\mtJht^{j,\theta} := \theta \mtJht +  (1-\theta ) \mtJht^{j-1}$ and the time derivatives are approximated by backward difference quotients.\\
The parameter $\theta$ is chosen as $\theta = \frac12 + \mathcaa{ O}(k)$, which corresponds to a shifted Crank-Nicolson scheme. By this choice one obtains second order accuracy in time and additionally recovers global stability \cite[Sec. 5.3]{RiWi}. The latter is important for stable 
behavior for long-term computations and not guaranteed by the standard Crank-Nicolson scheme, see \cite{Wi13}.

\subsection{Spatial Discretization}
\label{subsec::disch}

For the spatial discretization, we use a triangulation $\mathcaa T_{\sf h}$ of the domain $\sdom \Omega$ with 4262 vertices and 8225 cells $K$. For the sake of clarity, and since we focus on presenting the main ideas,
we denote the discretized domains also by $\tilde \Omega$, $\tilde \Omega_f$ and $\tilde \Omega_s$. Moreover, we also do not write the subscript ${\sf h}$ for discretized boundaries, e.g., we denote $\sdom \Gamma_{d, \sf h}$ by $\sdom \Gamma_d$. In order to have a stable discretization of the Navier-Stokes part of the FSI equations, we choose lowest order Taylor-Hood elements $ (\sdom \vel_{\sf h}, \sdom{p}_{\sf h} )\in  (\mathcaa P^2 (\mathcaa T_{\sf h})^d , \mathcaa P^1 (\mathcaa T_{\sf h} ) )$,
where
\[\mathcaa P^l (\mathcaa T_{\sf h})^m := \lbrace  \sdom \vel_{\sf h} \in \mathcaa C(\bigcup_{K \in \mathcaa T_{\sf h}} K)^m ~:~ \sdom \vel_{\sf h}\vert_K \text{ is a polynomial of degree }l, \ \forall K \in \mathcaa T_{\sf h} \rbrace \]
for $l \geq 0$ and $m \in \mathbb N$,
i.e., $\sdom \vel_{\sf h}$ is continuous and element-wise quadratic and $\sdom p_{\sf h}$ is continuous and linear on every element. Since $\sdom \vel_{\sf h}$ is equal to the temporal derivative of $\sdom \dis_{\sf h}$ on $\sdom \Omega_s$, $\sdom \dis_{\sf h}$ is chosen such that it has the same degrees of freedom as $\sdom \vel_{\sf h}$. Therefore, we choose $ (\sdom \dis_{\sf h}, \sdom \disa_{\sf h} ) \in  (\mathcaa P^2 (\mathcaa T_{\sf h})^d , \mathcaa P^2 (\mathcaa T_{\sf h})^d  )$. 
More precisely, the spaces $\sdom{\mathbf V}_\mmtau$, $\sdom{\mathbf V}_{0, \mmtau}$, $\sdom{\mathbf W}_\mmtau$, $\sdom{\mathbf W}_{0, \mmtau}$, $\sdom{\mathbf Z}_\mmtau$, and $\sdom P_\mmtau$ are approximated by the spaces
\begin{align*}
& \sdom{\mathbf V}_{\sf h} =  \lbrace \sdom \mv_{\sf h} \in \mathcaa P^2 ( \mathcaa T_{\sf_h})^d   ~:~  \sdom \mv\vert_{\sdom \Gamma_{fD}} = \sdom \mv_{fD,{\sf h}}  \rbrace, \\
& \sdom {\mathbf V}_{0, {\sf h}} = \lbrace \sdom \mv_{\sf h} \in \mathcaa P^2 ( \mathcaa T_{\sf_h})^d  ~:~ \sdom \mv_{\sf h}\vert_{\sdom \Gamma_{fD}} = 0  \rbrace, \\
& \sdom {\mathbf W}_{\sf h} =  \lbrace \sdom \dis_{\sf h} \in \mathcaa P^2 ( \mathcaa T_{\sf_h})^d  ~:~ \sdom \dis_{\sf h}\vert_{\sdom \Gamma_{fD} \cup \sdom \Gamma_{fN} } = 0, ~ \sdom \dis_{\sf h}\vert_{\sdom \Gamma_{sD}} = \sdom \dis_{sD, {\sf h}}  \rbrace, \\
& \sdom {\mathbf W}_{0, {\sf h}}  =  \lbrace \sdom \dis_{\sf h} \in \mathcaa P^2 ( \mathcaa T_{\sf_h})^d  ~:~ \sdom \dis_{\sf h}\vert_{\sdom \Gamma_{fD} \cup \sdom \Gamma_{fN}} = 0, ~ \sdom \dis_{\sf h}\vert_{\sdom \Gamma_{sD}} = 0  \rbrace, \\
& \sdom {\mathbf Z}_{\sf h}  =  \mathcaa P^2 ( \mathcaa T_{\sf_h})^d, \\
& \sdom P_{\sf h} = \mathcaa P^1 ( \mathcaa T_{\sf_h}).
\end{align*}

\subsection{Discretization of Objective Function}
\label{subsec::discof}

The spatial discretization of the objective function is determined by the discretization of the state of the FSI problem. To discretize the appearing time derivative terms, we use a finite difference scheme, more precisely, the time derivative $\partial_\mt \sdom \vel_{{\sf h}} (t_j)$ is approximated by $ (t_j - t_{j-1})^{-1} (\sdom \vel_{{\sf h}} (t_j) - \sdom \vel_{{\sf h}} (t_{j-1}) )$. The time integral is approximated using the trapezoidal rule.

\subsection{Discretization of Shape Transformations}
\label{subsec::discst}
In Section \ref{sec::shapetrans}, it is motivated that the choice of admissible shape transformations is delicate and requires available existence and regularity theory for the governing PDEs. However, existence and regularity theory for FSI systems is only available for special cases and under additional restrictions or assumptions.  In particular, there are no theoretical results concerning existence and regularity of solutions available for the model \eqref{coupledsystem3}. Thus, we restrict the considerations to the discretized problem. Here, the main requirements for choosing admissible shape transformations reduce to ensure the following:
\begin{itemize}
	\item The source term, the boundary and initial conditions remain untouched by admissible shape transformations. 
	\item The transformed triangulation $\mmtau_{\sf h}(\mathcaa T_{\sf h})$ is the discretization of a Lipschitz domain, which means that mesh degeneration is prevented. This is a delicate task that 
	gained attention in several publications. 
	In the context of shape optimization see, e.g., \cite{IgStWe} and the references therein, 
	in the context of ALE transformations see, e.g., \cite{B17, EF16, HK22}.
	Mesh degeneration is not seen directly since all computations are performed on the fixed shape reference $\sdom \Omega$ domain, however, it is the main bottleneck in the performance of the optimization.
    In this work, we tackle it by choosing the set of 
	admissible transformations carefully from a continuous perspective, see Section \ref{subsec::adtrans}.  
	\item The space $\sdom {\mathbf Y}_{\sf h} (\mathcaa T_{\sf h}) \circ \mmtau_h$ is isomorphic to $\sdom {\mathbf Y}_h (  \mmtau_h(\mathcaa T_{\sf h}))$ for all $\mmtau_{\sf h} \in \mTadh$, where $\sdom {\mathbf Y}_{\sf h}$ denotes the discrete state space and $ \mTadh$ the discrete set of admissible transformations. To do so, we choose $\mTadh \subset \mathcaa P^1 (\mathcaa T_{\sf h})$.
	\item There is a one-to-one-correspondence between transformations and shapes. Analogously to the continuous case, we choose a scalar valued variable $\tilde d_{\sf h} \in \sdom D_{\Gamma_d, {\sf h}}$, where $\sdom D_{\Gamma_d, {\sf h}}$ denotes the space of piecewise linear functions on $\Gamma_d$, in addition, require that $\mmtau_{\sf h}$ is equal to the identity on $\sdom \Omega_{s}$ and consider the discretized version of the operator $\mathrm B$ presented in Section \ref{subsec::adtrans}.
\end{itemize}

\subsection{Discretized Version of the Shape Optimization Problem}
\label{subsec::dop}
Let $\mathbf d \in \mathbb R^{n}$ denote the vector that represents $\sdom d_h$ via
\[
\sdom d_h=\Psi(\mathbf d) := \sum_{i = 1}^{n} {\mathbf d}_i \Phi_i,
\]
where $\{\Phi_i~:~i=1,\ldots,n\}$ is the nodal basis of the piecewise linear continuous elements on $\sdom D_{\Gamma_d}$.
Moreover, let $ \mathrm B_{\sf h}: \sdom D_{\Gamma_d, {\sf h}} \to \mathcaa P^1(\mathcaa T_{\sf h})^d$, $\sdom d_h \mapsto \sdom{ \mathbf u}_{\tau, {\sf h}}$, be the discretization of $\mathrm B$ which is obtained by solving the discretizations of \eqref{LapBel} and \eqref{extension::equation}, and $\sdom j_\Omega$ and $ \sdom g_\Omega$ be defined as in Section
\ref{subsec::op}.
The discretized shape optimization problem then is
\begin{align}
	\begin{split}
\min_{\mathbf d \in \mathbb R^{{n}} } & f(\mathbf d) + \frac{\alpha}2 \mathbf d^\top \mathbf S \mathbf d\quad
\mathrm{s.t. } \quad g(\mathbf d) = 0,
\end{split}
\label{generalform}
\end{align}
where
\begin{align}
 {\mathbf S}_{i,j}=
( \Phi_i, \Phi_j)_{\sdom D_{\Gamma_d, \sf h}} = (\Phi_i, \Phi_j )_{L^2(\sdom \Gamma_d)}.
\label{definition::S}
\end{align}
The vector $\mathbf d$ characterizes a transformation via the following chain of compositions
\begin{align*}
\mathbf d \xmapsto{~\Psi~} \sdom d_{\sf h} \xmapsto{~\mathrm B_{\sf h}~} \sdom{\mathbf u}_{\tau, {\sf h}} \xmapsto{~\mathrm{id}_\mz +\,\text{\bf\boldmath$\cdot$}~}
\mmtau_{\sf h}.
\end{align*}
The objective is defined by $f(\mathbf d):= \sdom j_\Omega (\sdom{\mathbf u}_{\tau, {\sf h}}) + \mathcal R_p(\sdom d_{\sf h})$ and $g(\mathbf d) = \sdom g_\Omega(\sdom{\mathbf u}_{\tau, {\sf h}})$, where $
\sdom{\mathbf u}_{\tau, {\sf h}} = \sdom{\mathbf u}_{\tau, {\sf h}}(\mathbf d)$, cf. \eqref{cont::opt::1}.

\section{Numerical Implementation}
\label{sec::fsinr}
The numerical tests\footnote{The source code is available at \url{https://github.com/JohannesHaubner/ShapeOpt/}.}  presented here are implemented in FEniCS \cite{ FenicsBook}, a collection of free software for the automated solution of PDEs.
For the computation of the gradients the additional package dolfin-adjoint \cite{dolfinadjoint} is used, which provides the automated differentiation of the reduced cost functional based on adjoint computations on the discrete system. For large scale computations, the simulation needs to be based on a checkpointing strategy \cite{GrWa00}, meaning that, to save memory, the forward solution is not saved for every time-step but only at several checkpoints. In order to solve the backwards equations, the forward equation is resolved starting from the closest checkpoints and then adjoint timesteps are made over this strip of recomputed states. The software package IPOPT \cite{ipopt} is used for solving the constrained optimization problem on the shape reference domain $\tilde\Omega$. In this paper, we solve the PDE system that we transformed by hand to $\tilde\Omega$. We mention, however, that there exist shape differentiation tools that automate this and build on the same transformation idea \cite{dokken2020automatic, ham2019automated}.

Many existing implementations of optimization methods, such as IPOPT, work with the Euclidean inner product. Therefore, handing the discretized optimization problem \eqref{generalform} directly to IPOPT leads to a loss of information since it is no longer taken into account that $\sdom d_{\sf h}$ is  the discretization of a function in a specific function space.
Since we have
$\sdom d_{\sf h}(\mz) := \sum_{i=1}^n  {\mathbf d}_i \Phi_i (\mz)$,
the correct inner product of vectors $\mathbf d^1$ and $\mathbf d^2$ is thus given by
$(\mathbf d^1)^{\top}  {\mathbf S} (\mathbf d^2)$,
where ${\mathbf  S }$ is defined in \eqref{definition::S}. 
Working on the space of transformed coordinates
$
\breve{\mathbf d} = \breve{\mathbf S} \mathbf d,
$
where $\breve {\mathbf S}$ is chosen such that $\breve{\mathbf S}^{\top} \breve{\mathbf S} = \mathbf S$, e.g., $\breve{\mathbf S} = \mathbf S^\frac12$ (which is impractical if the size of $\mathbf S$ is large) or obtained by a (sparse) Cholesky decomposition, takes the above considerations into account. 
We pass the following functions to IPOPT $$ \breve f:  \mathbb R^{n} \to \mathbb R, \quad \breve{\mathbf{ d}} \mapsto f ( \breve{\mathbf S}^{-1} \breve{\mathbf{ d}}) + \frac{\alpha}2 \breve{ \mathbf d}^{\top} \breve{ \mathbf d}, $$
as well as,
$ \nabla \breve f : \mathbb R^{n} \to \mathbb R^{n},$ $\breve{\mathbf{ d}} \mapsto \breve{\mathbf S}^{-\top} \nabla f ( \breve{\mathbf S}^{-1} \breve{\mathbf{ d}}) + \alpha \breve{\mathbf{ d}} .$
This has several advantages in the numerical solution process of the optimization problem. Performing a steepest descent method with step size $1$ for the function $\breve f$ results in $\mathbf d^k =\breve{ \mathbf S}^{-1} \breve{\mathbf{ d}}^k$ and
$$ {\mathbf d}^{k+1} = \breve {\mathbf S}^{-1} \breve{\mathbf{ d}}^{k+1} = \breve {\mathbf S}^{-1} ( \breve{\mathbf{ d}}^k - (\breve{\mathbf S}^{-\top} \nabla f ( \breve{\mathbf S}^{-1}  \breve{\mathbf{d}}^k) + \alpha \breve{\mathbf{ d}}^k))  = \mathbf d^k - (\mathbf S^{-1} \nabla f (\mathbf d^k) + \alpha \mathbf d^k),$$ i.e., a steepest descent method for $f(\cdot) + \frac\alpha{2} (\cdot)^\top \mathbf S (\cdot)$ using the Riesz representation of the derivative, which is given by $\mathbf S^{-1} \nabla f(\cdot) + \alpha (\cdot)$. In \cite{ScFuHa} it is shown that this leads to mesh-independent convergence rates for some examples. This is also expected for other optimization algorithms. Hence, we consider the optimization problem
\begin{align*}
\min_{ \breve{\mathbf d} \in \mathbb R^{{n}} } & f( \breve{\mathbf S}^{-1} \breve{\mathbf d}) + \frac{\alpha}2 \breve{ \mathbf d}^{\top} \breve{ \mathbf d} \quad
\mathrm{s.t. } \quad \breve g(\breve{\mathbf d}) = 0,
\end{align*}
where $\breve{g } (\breve{ \mathbf d}) = g ( \breve{ \mathbf S}^{-1} \breve{ \mathbf d})$. Alternatively, one could apply an algorithm which directly works with the correct inner product.

We remark that, in our case, $\sdom D_{\Gamma_d}$ is chosen as  $L^2(\sdom \Gamma_d)$, thus ${\mathbf S}$ is a mass matrix that is spectrally equivalent to the lumped mass matrix, which is diagonal. Therefore, our change of variables has a similar effect as suitably rescaling the variables. But if smoother or less smooth transformations would be desired, then $\sdom D_{\Gamma_d}$ could be chosen as some other space, e.g,
$H^s(\sdom \Gamma_d)$, with $s>0$ (or $<0$) corresponding to more (or less) smoothness.

\section{Numerical Results on Basis of the FSI2 Benchmark}
\label{sec::numres}

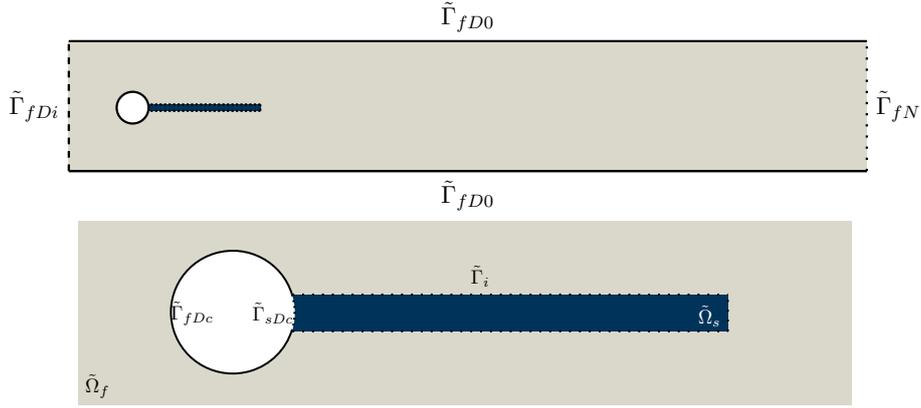
\begin{figure}
	\centering
	\scalebox{.84}{
	\begin{tikzpicture}
	\fill[tumg] (0.0,0.0) rectangle (5*2.5, 5*0.41);
	\fill[tumblues2] (5*0.2, 5*0.19) rectangle (5*0.6, 5*0.21);
	\draw[densely dotted, line width = 1pt] (5*0.2, 5*0.19) -- (5*0.6, 5*0.19);
	\draw[densely dotted, line width = 1pt] (5*0.2, 5*0.21) -- (5*0.6, 5*0.21);
	\draw[densely dotted, line width = 1pt] (5*0.6, 5*0.19) -- (5*0.6, 5*0.21);
	\fill[white] (5*0.2,5*0.2) circle (5*0.05);
	\draw[line width = 1pt] (0,0) -- (5*2.5,0);
	\draw[line width = 1pt] (0,0.41*5) -- (5*2.5,0.41*5) node[midway, above] {$\sdom  \Gamma_{fD0}$};
	\draw[line width = 1pt] (0,0) -- (5*2.5,0) node[midway, below] {$\sdom  \Gamma_{fD0}$};
	\draw[loosely dotted, line width = 1pt] (2.5*5, 0.0) -- (2.5*5, 5*0.41) node[midway, right] {$\sdom  \Gamma_{fN}$};
	\draw[dashed, line width = 1pt] (0.0,0.0) -- (0.0,5*0.41) node[midway, left] {$\sdom  \Gamma_{fDi}$};
	\draw[line width = 1pt] (5*0.2, 5*0.2) circle (5*0.05);
	\end{tikzpicture}
	}
    \centering
    \scalebox{.74}{
	\begin{tikzpicture}
	\small
	\centering
	\fill[tumg] (0.0,0.0) rectangle (22*0.625,22*0.15);
	\fill[tumblues2] (22*0.125, 22*0.060) rectangle (22*0.525, 22*0.090);
	\draw[loosely dotted, line width = 1pt] (22*0.125, 22*0.06) -- (22*0.525, 22*0.06);
	\draw[loosely dotted, line width = 1pt] (22*0.125, 22*0.09) -- (22*0.525, 22*0.09) node[midway,above]{$\sdom  \Gamma_i$};
	\draw[loosely dotted, line width = 1pt] (22*0.525, 22*0.06) -- (22*0.525, 22*0.09);
	\fill[white] (22*0.125, 22*0.075) circle (22*0.05);
	\draw[line width = 1pt] (22*0.172697, 22*0.09) arc(atan(0.015/0.05):(360-atan(0.015/0.05)):22*0.05);
	\node[right] at (0.135*22,0.075*22) {$\sdom  \Gamma_{sDc}$};
	\draw[line width = 1pt, dotted] (22*0.172697, 22*0.09) arc((atan(0.015/0.05)):(-atan(0.015/0.05)):22*0.05);
	\node[left] at (0.115*22,0.075*22) {$\sdom  \Gamma_{fDc}$} ;
	\node[above right] at (0.0,0.0) {$\sdom \Omega_f$};
	\node[above left, color = white] at (22*0.525,22*0.058) {$\sdom \Omega_s$};
	\end{tikzpicture}
	}
	\caption{Shape reference domain $\sdom \Omega$}
	\label{fig::fsi::shaperef}
\end{figure}

We build the validation of our numerical implementation on the FSI2 benchmark, which was proposed in \cite{TuHr}. This benchmark considers the coupling of the Navier-Stokes 
equations and Saint Venant-Kirchhoff type material equations in a two-dimensional rectangular domain of length $l = 2.5$ and height $h = 0.41$, the bottom left corner of which is located at the origin $(0,0)^{\top}$. On the left boundary $\aledom \Gamma_{fDi}$ we have a parabolic inflow given by 
\begin{align*}
\sdom {\mathbf v}_{fD} ((0,\mz_2)^{\top}, \mt ) = \begin{cases}
(3 \bar v h^{-2} \mz_2 (h-\mz_2 )(1 - \mathrm{cos} (\frac{\pi}{2} \mt )),0)^{\top} \quad & \text{if } \mt< 2.0, \\
(6 \bar v h^{-2} \mz_2 (h-\mz_2 ),0)^{\top} &\text{otherwise},
\end{cases}
\end{align*}
with mean inflow velocity $\bar v$. No-slip conditions on the bottom and top $\sdom \Gamma_{fD0}$ and do-nothing boundary conditions on the right boundary $\sdom \Gamma_{fN}$, i.e., $\sdom{\mathbf g}_f = 0$, are imposed. In this pipe, there is a circular obstacle with radius $r = 0.05$ centered at $ (0.2,0.2 )$ to which an elastic beam of length $0.4$ and width $0.02$ is attached as illustrated in Figure \ref{fig::fsi::shaperef}. 
On $\sdom \Gamma_{fDc}$, the part of the circular obstacle that is part of the fluid boundary, homogeneous Dirichlet boundary conditions are imposed on the fluid velocity and solid displacement. The initial conditions are set to $0$. Thus, the fluid-structure system is completely determined by the parameters
$\rho_s = 1\cdot 10^{4}$, $\lambda_s = 2 \cdot 10^6$, $\mu_s = 5 \cdot 10^5$, $\rho_f = 1\cdot 10^3$, $\nu_f = 1 \cdot 10^{-3}$, $\bar v = 1$ and $\alpha_p = \alpha_w = 1 \cdot 10^{-9}$. Details on the implementation are given in the previous sections. As time horizon for the optimization of the mean drag we choose $T = 15 s$. Given a design boundary (or interface) $\sdom \Gamma_d$, we want to optimize the shape of $\sdom \Omega$ such that the fluid drag is minimized.

\subsection{Optimization of Shape of Obstacle}
\label{sec::61}

\begin{figure}
	\captionsetup[subfloat]{farskip=2pt,captionskip=1pt}
	\begin{tabular}{*{2}{b{0.5\textwidth-2\tabcolsep}}}
		\includegraphics[width=0.85\hsize]{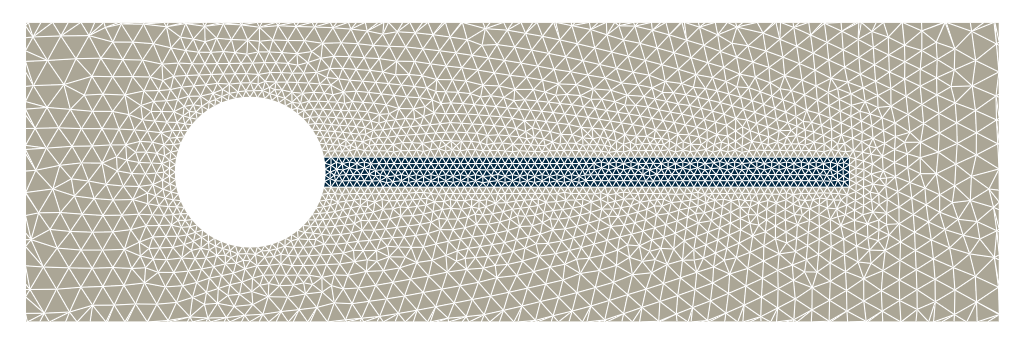}
		\subfloat
		{\includegraphics[width=0.85\hsize]{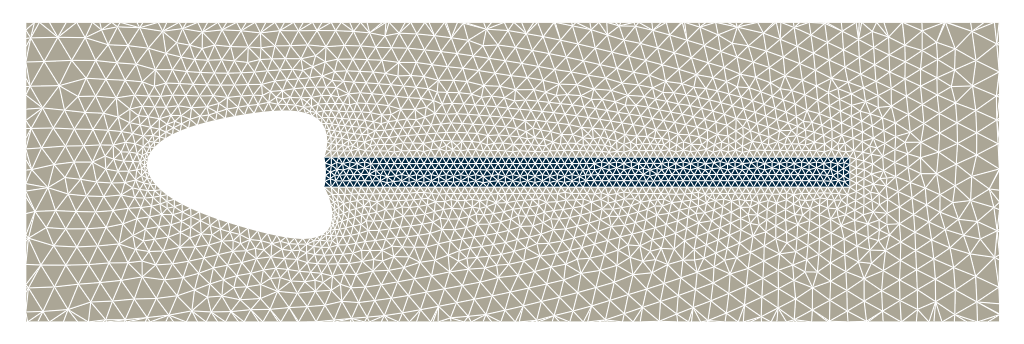}
		}\vspace{0.5cm}
		&
		\subfloat{
			\includegraphics[width=0.98\hsize]{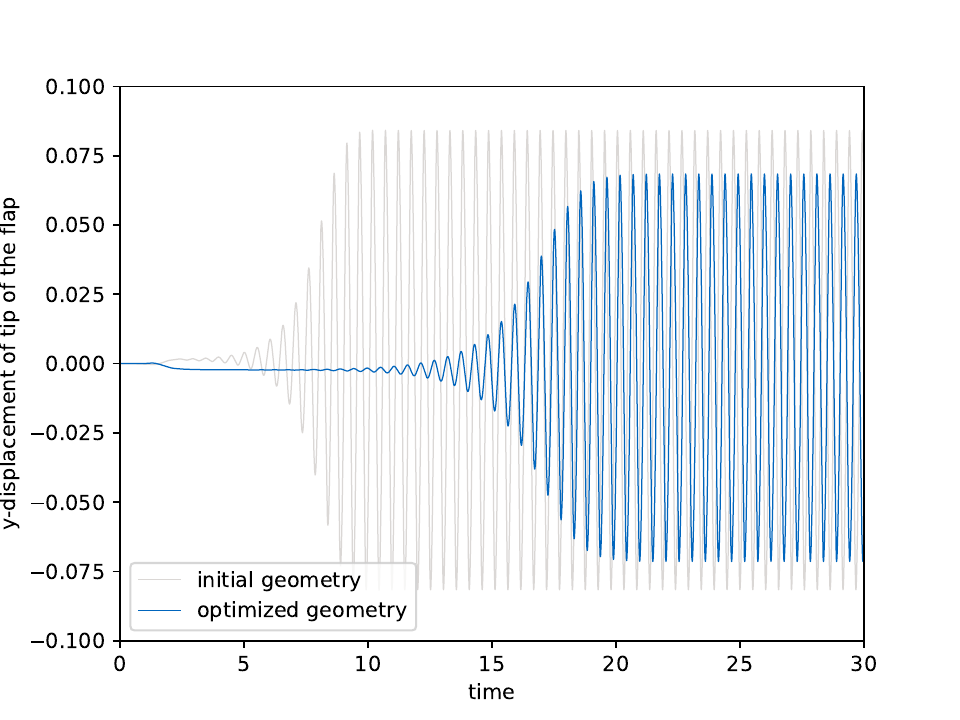}}
	\end{tabular}
	\caption{Optimization of shape of obstacle (Section \ref{sec::61}): Comparison of vertical displacement of the tip of the flap for the initial (top left) and optimized (bottom left) design}
	\label{fig::cio}
\end{figure}

For a first example, $\sdom \Gamma_d = \sdom \Gamma_{fDc}$ serves as design boundary, i.e., we optimize the shape of the circular obstacle and keep the solid domain fixed.
In addition, we use the regularization parameter $\alpha = 0.1$ and we choose $\beta=0.01$ in the equation \eqref{tildebetaeq} for the weighting $\tilde\beta$ of the Laplace-Beltrami operator in \eqref{LapBel}. As constraints, we require that the volume and barycenter of the obstacle remain the same as for the inital configuration. For our example, IPOPT converges after $14$ iterations with an overall scaled NLP error (cf. \cite[p. 3, (5)]{ipopt}) smaller than $1\cdot 10^{-3}$, see Table \ref{tab::or2}. The objective function value is reduced by more than $30$ \%. Figure \ref{fig::cio} compares the initial configuration and the optimized configuration. The vertical displacement of the tip of the flap shows that even though the optimization is only performed on the first $15$ seconds the amplitude of the vertical displacement of the tip of the flap is also smaller for long-term simulations. 

\begin{table}[ht!]
  \begin{tabularx}{\textwidth}{R|x|x|x|xe}
    \arrayrulecolor{white}
    \rowcolor{tumblues2}
    \textcolor{white}{iteration} & \textcolor{white}{objective}& \textcolor{white}{objective \mbox{ \small w/o reg. \& pen.}}& \textcolor{white}{dual infeasibility} & \textcolor{white}{linesearch-steps} &\\[0.5ex]
0 & $1.5985e+02$ & $1.5985e+02$ & $2.27e+01$ & $0$ & \\[0.5ex]
    \rowcolor{tumg}
1 & $1.1851e+02$ & $1.1569e+02$ & $9.62e+00$ & $2$ & \\[0.5ex]
2 & $1.1436e+02$ & $1.1060e+02$ & $4.40e+00$ & $1$ & \\[0.5ex]
    \rowcolor{tumg}
3 & $1.1302e+02$ & $1.0971e+02$ & $2.43e+00$ & $1$ & \\[0.5ex]
4 & $1.1162e+02$ & $1.0818e+02$ & $1.21e+01$ & $3$ & \\[0.5ex]
    \rowcolor{tumg}
5 & $1.0773e+02$ & $1.0393e+02$ & $5.75e+00$ & $1$ & \\[0.5ex]
6 & $1.0756e+02$ & $1.0349e+02$ & $5.06e+00$ & $1$ & \\[0.5ex]
    \rowcolor{tumg}
7 & $1.0751e+02$ & $1.0337e+02$ & $1.01e+00$ & $1$ & \\[0.5ex]
8 & $1.0747e+02$ & $1.0336e+02$ & $3.29e-01$ & $1$ & \\[0.5ex]
    \rowcolor{tumg}
9 & $1.0743e+02$ & $1.0330e+02$ & $3.69e-01$ & $1$ & \\[0.5ex]
10 & $1.0730e+02$ & $1.0316e+02$ & $6.93e-01$ & $1$ & \\[0.5ex]
    \rowcolor{tumg}
11 & $1.0726e+02$ & $1.0331e+02$ & $4.60e-01$ & $1$ & \\[0.5ex]
12 & $1.0725e+02$ & $1.0333e+02$ & $8.01e-01$ & $2$ & \\[0.5ex]
    \rowcolor{tumg}
13 & $1.0725e+02$ & $1.0332e+02$ & $8.10e-01$ & $2$ & \\[0.5ex]
14 & $1.0725e+02$ & $1.0329e+02$ & $4.33e-03$ & $1$ & \\[0.5ex]
    \rowcolor{tumg}
  \end{tabularx}
  \caption{Optimization of shape of obstacle (Section \ref{sec::61}): Optimization results when IPOPT converges up to an overall scaled NLP tolerance of $10^{-3}$. The table shows the objective function value, the objective function value without the regularization and penalization term, the dual infeasibility and the number of linesearch-steps per iteration.}
  \label{tab::or2}
\end{table}

\subsection{Optimization of Interface}
\label{sec::62}
\begin{figure}[h!]
	\captionsetup[subfloat]{farskip=2pt,captionskip=1pt}
	\begin{tabular}{*{2}{b{0.5\textwidth-2\tabcolsep}}}
        \cellcolor{white}
		\includegraphics[width=0.85\hsize]{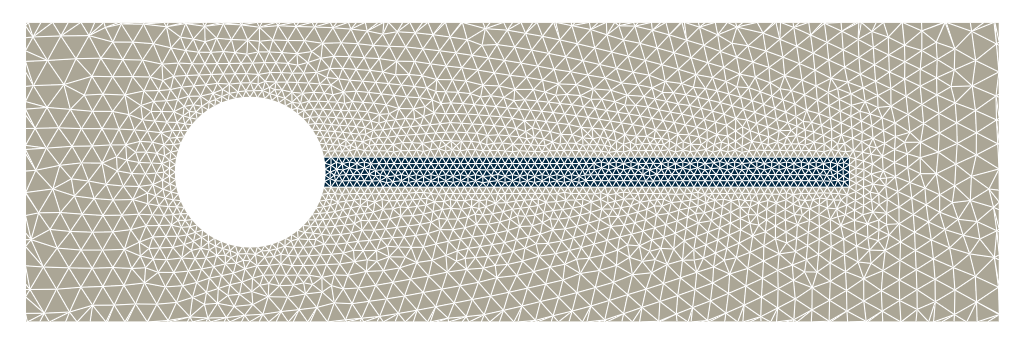}
		\subfloat
		{\includegraphics[width=0.85\hsize]{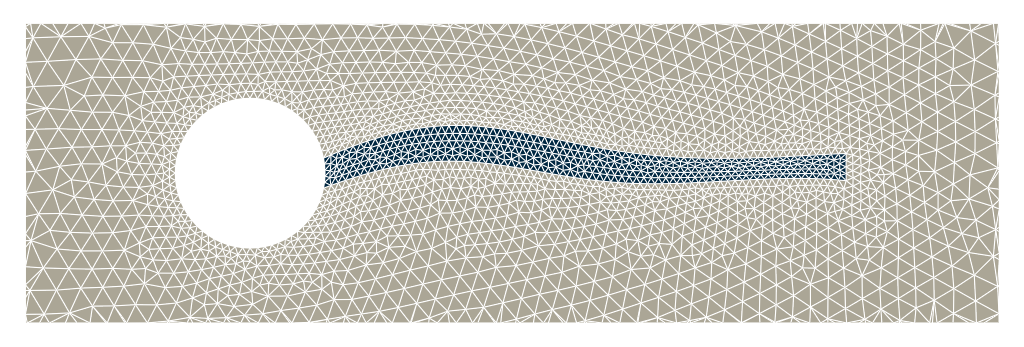}
		}\vspace{0.5cm}
		&
		\subfloat{
			\includegraphics[width=0.98\hsize]{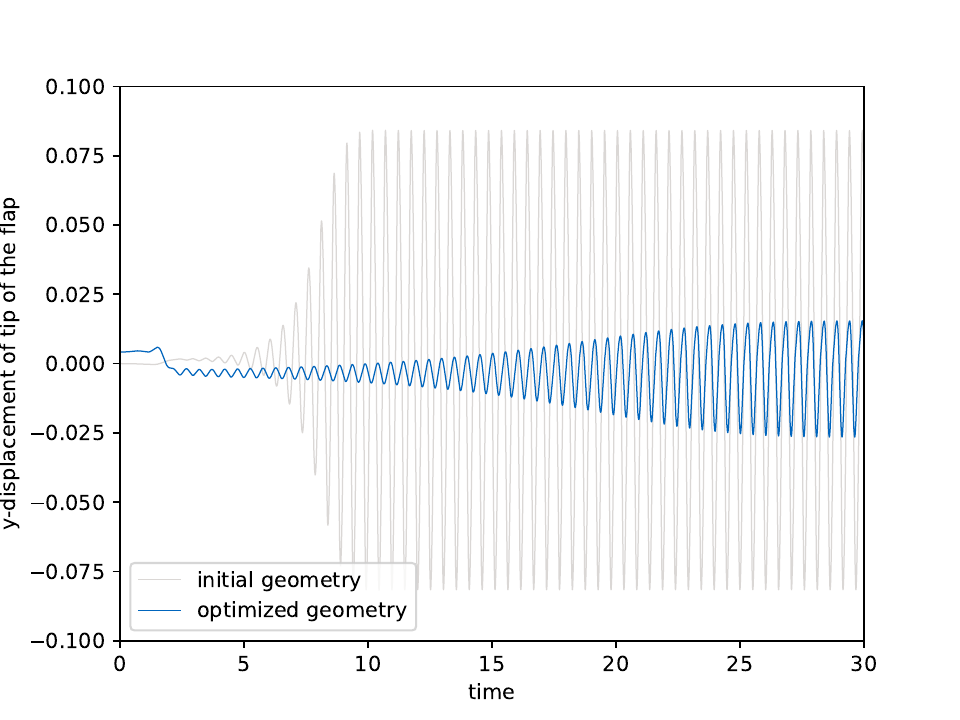}}
	\end{tabular}
	\caption{Optimization of interface (Section \ref{sec::62}): Comparison of vertical displacement of the tip of the flap for the initial (top left) and optimized (bottom left) design. The initial displacement of the tip of the flap for the optimized design is not plotted as zero, because the tip of the optimized flap is located at (0.597718, 0.20417) instead of (0.6,0.2). }
	\label{fig::cio2}
\end{figure}

For a second example, we choose the interface $\sdom \Gamma_d = \sdom \Gamma_{i}$ as design boundary.
Here, we use the regularization parameter $\alpha = 10$ and weighting for the Laplace-Beltrami operator $ \beta = 10^{-3}$. Too small values for the regularization parameter $\alpha$ lead to poor convergence. As constraint, we impose that the volume of the solid domain remains the same as for the initial configuration. For this example, IPOPT converges after $14$ iterations with an overall scaled NLP error smaller than $1\cdot 10^{-3}$, see Table \ref{tab::or22}. The objective function value is reduced by more than $20$ \%. As for the previous example, Figure \ref{fig::cio2} compares the initial configuration and the optimized configuration.

\begin{table}[ht!]
  \begin{tabularx}{\textwidth}{R|x|x|x|xe}
    \arrayrulecolor{white}
    \rowcolor{tumblues2}
    \textcolor{white}{iteration} & \textcolor{white}{objective}& \textcolor{white}{objective \mbox{ \small w/o reg. \& pen.}}& \textcolor{white}{dual infeasibility} & \textcolor{white}{linesearch-steps} &\\[0.5ex]
0 & $1.5985e+02$ & $1.5985e+02$ & $3.85e+01$ & $0$ & \\[0.5ex]
    \rowcolor{tumg}
1 & $1.4593e+02$ & $1.4479e+02$ & $1.67e+01$ & $10$ & \\[0.5ex]
2 & $1.3925e+02$ & $1.3592e+02$ & $1.46e+01$ & $1$ & \\[0.5ex]
    \rowcolor{tumg}
3 & $1.3246e+02$ & $1.2996e+02$ & $5.85e+00$ & $1$ & \\[0.5ex]
4 & $1.2932e+02$ & $1.2732e+02$ & $2.91e+00$ & $1$ & \\[0.5ex]
    \rowcolor{tumg}
5 & $1.2782e+02$ & $1.2612e+02$ & $1.37e+00$ & $1$ & \\[0.5ex]
6 & $1.2757e+02$ & $1.2567e+02$ & $3.63e-01$ & $1$ & \\[0.5ex]
    \rowcolor{tumg}
7 & $1.2754e+02$ & $1.2565e+02$ & $3.08e-01$ & $1$ & \\[0.5ex]
8 & $1.2744e+02$ & $1.2541e+02$ & $2.26e-01$ & $1$ & \\[0.5ex]
    \rowcolor{tumg}
9 & $1.2743e+02$ & $1.2528e+02$ & $3.08e-02$ & $1$ & \\[0.5ex]
10 & $1.2743e+02$ & $1.2528e+02$ & $1.99e-02$ & $1$ & \\[0.5ex]
    \rowcolor{tumg}
11 & $1.2743e+02$ & $1.2528e+02$ & $4.49e-03$ & $1$ & \\[0.5ex]
12 & $1.2743e+02$ & $1.2528e+02$ & $5.56e-03$ & $2$ & \\[0.5ex]
    \rowcolor{tumg}
13 & $1.2743e+02$ & $1.2528e+02$ & $4.85e-03$ & $1$ & \\[0.5ex]
14 & $1.2743e+02$ & $1.2527e+02$ & $9.04e-04$ & $1$ & \\[0.5ex]
    \rowcolor{tumg}
  \end{tabularx}
  \caption{Optimization of interface (Section \ref{sec::62}): Optimization results when IPOPT converges up to an overall scaled NLP tolerance of $10^{-3}$. The table shows the objective function value, the objective function value without the regularization and penalization term, the dual infeasibility and the number of linesearch-steps per iteration.}
  \label{tab::or22}
\end{table}

\begin{remark}
    Note that the optimal solutions depend on several choices such as, e.g., the parameterization of the set of admissible shapes via the boundary and extension operators, the choice of the regularization parameter and the threshold for the determinant constraint. The presented method can be used as a building block of an outer iterative procedure, where the (potentially remeshed) optimized geometry of the previous step is used as initial domain for the next step and every outer iterate fulfills the constraints up to a certain, user-specified tolerance. This can lead to refined results, see, e.g., \cite{HSU20} for Stokes flow. 
\end{remark}

\section{Conclusion and Outlook}

We applied the method of mappings to solve a shape optimization problem for unsteady FSI. We introduced the FSI model and show in the appendix that the adjoint equations for a linearized version of the model attain the same structure as the forward equations. Performing shape optimization via the method of mappings requires the transformation of these equations to a reference domain. The choice of admissible shape transformations is done sophisticated in order to have a well-posed optimization problem. One of the main challenges is the prevention of mesh degeneration during the optimization process, which is done with a formulation of the optimization method respecting the continuous requirements for the transformations and adding a penalization of too small determinant values of the gradient of the transformation. As objective function a volume formulation of the drag of the obstacle is chosen. 
Numerical results show the viability of the approach to optimize the shape of the obstacle's boundary and of the interface based on problems derived from the FSI2 benchmark. There are several degrees of freedom in realizing the method of mappings, e.g., via the choice of the boundary operator and the extension operator. Moreover, also iterative approaches can be considered, where the optimized shape serves as initial shape for the next iteration. This can also be combined with remeshing routines and parameter tuning. In addition, to make the algorithm more robust, the determinant constraint can also be imposed on an hold all domain that includes the obstacle. 

\section*{Acknowledgments}
This work was funded by the Deutsche
Forschungsgemeinschaft (DFG, German Research Foundation) 
as part of the International Research Training Group IGDK 1754
``Optimization and Numerical Analysis for Partial Differential Equations with Nonsmooth
Structures''--Project Number 188264188/GRK1754. Johannes Haubner also acknowledges support from the Research Council of Norway, grant 300305.
The authors would like to acknowledge J{\o}rgen S. Dokken, Henrik N. Finsberg and Simon W. Funke for the support and discussions on reproducibility of numerical results.

\appendix

\section{Adjoint for a Linear Unsteady FSI Problem with Stationary Interface}
\label{sec::adjcon}

For efficient derivative computations in optimal control, the adjoint equations have to be solved. Especially in cases where no automatic differentiation can be applied, it is crucial to derive an explicit formula for the adjoint equations. Even though the FSI model is modified for performing shape optimization, the adjoint equations to the unmodified model can be used to compute the derivative. This can be realized by performing every iteration on the current ALE reference domain instead of the nominal domain, cf. \cite[Sec. 2.2.2]{BaLiUl}. In our choice of the software framework, automatic differentiation is available and the following considerations are not needed for the numerical implementation. Nevertheless, we think that the result might be of interest for some readers.

We consider the adjoint of a linear version of the fluid-structure interaction model \eqref{coupledsystem}. More precisely, we consider Stokes flow for the fluid and linear elasticity for the solid equations. Additionally, we restrict ourselves to the case with a stationary interface $\aledom \Gamma_i$ and homogeneous Dirichlet boundary conditions, i.e., $\aledom \Omega = \pdom\Omega (\mt) = \Omega$ and $\partial \Omega(\mt) = \overline{\aledom \Gamma_{fD}} \cup \overline{\aledom \Gamma_{sD}}$ for any $\mt \in I$. This also implies that $\aledom {\mX} = \mathrm{id}_\my$,  ${\maledomJ} = 1$ and $\maledomJX = \mId$. The resulting fluid-structure interaction problem (for better readability without subscripts for the state variables and without superscripts) reads as follows
\begin{align}
\begin{split}
\mrhof \partial_\mt \vel - \mathrm{div}  (\msigmaf  ) = \mrhof \mathbf f_f  &\quad \text{in } \mQf, \\
\mathrm{div}  ( \vel )  = 0 &\quad \text{in } \mQf, \\
\mrhos \partial_\mt \vel - \mathrm{div} ( \msigmas  )  = \mrhos \mathbf f_s &\quad \text{in } \mQs,\\
\mrhos (\partial_\mt \dis - \vel )  = 0 &\quad\text{in } \mQs, 
\end{split}
\label{eq::linfsi::adj}
\end{align}
with the initial conditions
\begin{align*}
\vel \mevalo = \vel_{0} \quad \text{on } \Omega_f,\quad\quad \dis \mevalo = \dis_0 \quad \text{on } \Omega_s, \quad\quad
\vel \mevalo = \vel_0 \quad \text{on } \Omega_s,
\end{align*}
boundary conditions
\begin{align*}
\vel = 0 \quad \text{on } \mGf, \quad\quad \dis = 0 \quad \text{on } \mGs, 
\end{align*}
and the additional coupling conditions
\begin{align*}
\partial_\mt \dis = \vel   \quad \text{on } \mGi, \quad \quad
- \msigmaf \mathbf n_f = \msigmas \mathbf n_s  \quad\text{on } \mGi,
\end{align*}
where $\msigmaf =  \mu_f  (D \vel + {D \vel}^{\top} ) - p \mId$, $\mu_f = \rho_f \nu_f$, $
\msigmas = \mu_s  (D \dis + {D \dis}^{\top} ) +  \lambda_s \mathrm{div} ( \dis )\mId $ and, for the sake of convenience, we introduced $\vel_0$ defined by $\vel_{0} \vert_{\Omega_f} = \vel_{0f}$ and $\vel_0 \vert_{\Omega_s} = \dis_1$. For compatibility reasons there holds $\mathbf w_0 \vert_{\Gamma_s} = 0$. This corresponds to the setting considered in \cite{DuGu}.

We are interested in the structure of the adjoint equations and therefore do calculations on a formal level in order to derive a formulation for the adjoint system. In particular, we do not analyze the regularity of solutions but only assume that all functions are smooth enough such that the appearing terms and operations are well-defined. For the analysis of \eqref{eq::linfsi::adj} we refer to \cite{DuGu, DuGu2, FaMeVe}. 
As usual for unsteady problems, the flow of information in the adjoint equation is reversed in time. In order to transfer the theory and numerical methods from the forward model to the adjoint, it is desirable that, except for time reversal, the adjoint system has a similar or the same structure as the forward model. In our approach the adjoint has the same structure as the FSI system \eqref{eq::linfsi::adj}.
In \cite{FaMeVe}, it is shown that a straightforward weak formulation can fail to have the desired property, basically due to the equation $\partial_{ \mt } \dis - \vel = 0$. 
As a remedy, it is proposed in \cite{FaMeVe} to work with $\nabla \partial_{ \mt } \dis - \nabla \vel = 0$ instead.
In the following, we apply ideas from \cite{DuGu} to reformulate the weak formulation and obtain an analogous result.\\
We introduce $\mU  \mevalt = \mathbf w_0 + \int_0^\mt \vel \mevals  d  \ms$. Since $\dis\mevalo= \mathbf w_0$ and $\partial_{ \mt } \dis - \vel = 0$ on 
$\Omega_s$, we can substitute $\dis \mevalt =  \dis_0 + \int_0^\mt \vel \mevals  d  \ms = \mU  \mevalt$ on $\Omega_s$ into \eqref{eq::linfsi::adj}.
If $\dis$ is smooth enough such that the trace $\partial_{\mt} \dis \vert_{\mGs}$ exists, we obtain from $\partial_\mt \dis - \vel = 0$ on $\mQs$ and $\dis\vert_{\mGs} = 0$ that $\vel\vert_{\mGs} = 0$.  
One can check that
$\mrhos(\partial_\mt \mU - \vel) = 0$ on $\mQs$, $\mU = 0$ on $\mGs$, and $\mU \mevalo = \mathbf w_0$ are satisfied by the definition of $\mU$. Moreover, $\partial_\mt \mU = \vel$ on $\mGi$ is satisfied if we require $\vel\mevalt \in H^1_0 (\Omega )$ for almost all $\mt \in I$, which implies uniqueness of the trace. 
Thus, system \eqref{eq::linfsi::adj} with its initial, boundary and coupling conditions is equivalent to:
\begin{align}
\begin{split}
\mrhof \partial_\mt \vel - \mu_f \mathrm{div}   (D \vel +{D \vel}^{\top}  ) + \nabla p = \mrhof \mathbf f_f  & \quad \text{in }\mQf, \\
\mathrm{div}  ( \vel  ) = 0 & \quad  \text{in } \mQf, \\
\vel = 0 & \quad \text{on } \mGf, \\
\mrhos \partial_\mt \vel - \mu_s \mathrm{div}   (D(\mU) + {D(\mU)}^{\top} ) - \lambda_s \nabla  (\mathrm{div}  ( \mU ) ) = \mrhos\mathbf f_s & \quad\text{in }\mQs,\\
\vel   = 0 & \quad \text{on } \mGs, \\
\vel \mevalo  = \mathbf v_0 &  \quad\text{on } \Omega,
\end{split}
\label{system::lin}
\end{align}
with the additional coupling condition 
\begin{align*}
p \mnormalf - \mu_f  (D \vel + {D \vel}^{\top}  ) \mnormalf &= \mu_s (D(\mU) +{D(\mU)}^{\top} ) \mnormals + \lambda_s \mathrm{div}  ( \mU  ) \mnormals ~  \text{on } \mGi .
\end{align*}
The following notation is used:
\begin{itemize}
	\item $(p, q)_{\Omega} := \int_{ \Omega } p q d \xi $ for all $p,q \in L^2(\Omega)$, $(\vel, \vela)_{\Omega} := \int_{\Omega} \vel \cdot \vela d \xi$ for all $\vel, \vela \in L^2(\Omega)^d$ and $(\mathbf A, \mathbf B)_\Omega := \int_\Omega \mathbf A: \mathbf B d \xi$ for all $\mathbf A, \mathbf B \in L^2(\Omega)^{d \times d}.$
	\item $(\vel, \vela)_{\Gamma} := \int_{ \Gamma } \vel \cdot \vela d S(\xi)$ for all $\vel, \vela \in L^2(\Gamma)^d$, where $ dS(\xi)$ is the surface measure on $\Gamma$.
	\item $\mdl p,q \mdr_{Q^T}:= \int_0^T (p(\cdot, \mt),q(\cdot, \mt))_{\Omega} d\mt$ for all $p,q \in L^2((0,T), L^2(\Omega))$, \\$\mdl \vel, \vela \mdr_{Q^T} := \int_{0}^T ( \vel \mevalt , \vela \mevalt )_ \Omega d\mt$ for all $\vel, \vela \in L^2((0,T),L^2(\Omega)^d)$ and \\$\mdl \mathbf A, \mathbf B\mdr_{Q^T} := \int_0^T (\mathbf A, \mathbf B)_ \Omega d \mt$ for all $\mathbf A, \mathbf B \in L^2 ((0,T), L^2(\Omega)^{d \times d}).$	
	\item $
	a_f (\vel,\disa ) = \frac{\mu_f}{2}  ( D \vel + {D \vel}^{\top}, D \disa + {D \disa}^{\top}  )_{\Omega_f}$, \\
	$a_s (\vel,\disa ) = \frac{\mu_s}{2}  ( D \vel + {D \vel}^{\top}, D \disa + {D \disa}^{\top} )_{\Omega_s} + \lambda_s  ( \mathrm{div}  (\vel ), \mathrm{div} ( \disa  )  )_{\Omega_s},
	$
	for $\vel, \disa \in H_0^1(\Omega)^d$.
\end{itemize}
For vector valued variables $\vel$ and $\disa$, we use the identities
\begin{align}
    \begin{split}
&(\mathbf A, D \disa)_\Omega = - (\mathrm{div}(\mathbf A), \disa)_\Omega
+ (\mathbf A \mathbf n, \disa)_{\partial\Omega}, \\
&(\mathrm{div} (\vel), \mathrm{div} (\disa))_{\Omega} = 
-(\nabla \mathrm{div}(\vel), \disa)_{\Omega} 
+ (\mathrm{div}(\vel) \mathbf n, \disa)_{\partial \Omega}.
\end{split}
\label{eq::id}
\end{align}
Testing \eqref{system::lin} with $(\mpsiv, \mpsip)$ such that $ \mpsiv \vert_{\partial \Omega \times I} = 0$ gives:
\begin{align}
	\begin{split}
&\mrhof  \mdl \partial_\mt \vel, \mpsiv  \mdr_\mQf- \mu_f  \mdl \mathrm{div}  (D \vel + {D \vel}^{\top} ), \mpsiv  \mdr_\mQf\\& +  \mdl\nabla p, \mpsiv \mdr_\mQf-  \mdl \rho_f \mathbf f_f, \mpsiv \mdr_\mQf+ \rho_s  \mdl \partial_\mt \vel, \mpsiv  \mdr_\mQs  \\
& -  (( \rho_s \mathbf f_s, \mpsiv  ))_\mQs  -  \mdl \mathrm{div} ( \vel ), \mpsip  \mdr_\mQf \\ &+ \mrhof ( \vel\mevalo - \vel_0, \mpsiv\mevalo  )_{\Omega_f} + \mrhos ( \vel\mevalo - \vel_0, \mpsiv\mevalo  )_{\Omega_s} \\
&  - \mu_s  \mdl \mathrm{div}  (D(\mU) + {D(\mU)}^{\top}  ) , \mpsiv  \mdr_\mQs-\lambda_s  \mdl \nabla  (\mathrm{div}  ( \mU ) ), \mpsiv  \mdr_\mQs = 0.
\end{split}
\label{weakformeq}
\end{align}
With \eqref{eq::id}, $\vel \vert_{\partial \Omega \times I} = 0$, and the identity $(\mathbf A^\top, \mathbf B^\top)_\Omega = (\mathbf A, \mathbf B)_\Omega$ we obtain the formulas
\begin{align*}
&  \mu_f  ( \mathrm{div}  (D \vel + {D \vel}^{\top}  ), \mpsiv  )_{\Omega_f} = - a_f (\vel,\mpsiv ) + \mu_f ( (D \vel + {D \vel}^{\top}  ) \mnormalf , \mpsiv )_{\Gamma_i},\\
&( \nabla p, \mpsiv  )_{\Omega_f} = -  ( p, \mathrm{div} (\mpsiv ) )_{\Omega_f} +  ( p \mnormalf, \mpsiv )_{\partial \Omega_f} = -  ( p, \mathrm{div} (\mpsiv ) )_{\Omega_f} +  ( p \mnormalf, \mpsiv )_{\Gamma_i}, 
\end{align*}
and 
\begin{align*}
	& \mu_s  (  \mathrm{div}  (D\dis + {D \dis}^{\top}  ), \mpsiv )_{\Omega_s} + \lambda_s  ( \nabla ( \mathrm{div} ( \dis) ), \mpsiv  )_{\Omega_s} \\
	& = - \frac{\mu_s}{2}  ( D \dis + {D \dis}^{\top} , D \mpsiv + {D \mpsiv}^{\top}  )_{\Omega_s} + \mu_s  ( (D \dis + {D \dis}^{\top}  ) \mnormals , \mpsiv )_{\Gamma_i}  \\
	& \hspace{12pt} - \lambda_s (  \mathrm{div} ( \dis) , \mathrm{div} (\mpsiv)  )_{\Omega_s} + \lambda_s (  \mathrm{div} ( \dis) \mnormals , \mpsiv  )_{\Gamma_i} \\
	& = - a_s (\dis,\mpsiv )  + \mu_s  ( (D \dis + {D \dis}^{\top}  ) \mnormals , \mpsiv  )_{\Gamma_i} + \lambda_s (  \mathrm{div} ( \dis) \mnormals , \mpsiv  )_{\Gamma_i}.
\end{align*}
Thus, \eqref{weakformeq} can be reformulated as 
\begin{align*}
&\rho_f  \mdl \partial_\mt \vel, \mpsiv  \mdr_\mQf+ \int_0^T a_f (\vel, \mpsiv  ) d  \mt - \mdl p, \mathrm{div}( \mpsiv)  \mdr_\mQf\\
&+\rho_s \mdl \partial_\mt \vel, \mpsiv  \mdr_\mQs + \int_0^T a_s (\mU, \mpsiv ) d  \mt -  \mdl \mathrm{div}( \vel), \mpsip \mdr_\mQf\\
&- \mdl \rho_f \mathbf f_f, \mpsiv  \mdr_\mQf-  \mdl \rho_s \mathbf f_s, \mpsiv  \mdr_\mQs \\ &+ \mrhof ( \vel \mevalo - \vel_0, \mpsiv \mevalo  )_{\Omega_f}  + \mrhos ( \vel \mevalo - \vel_0, \mpsiv \mevalo  )_{\Omega_s} \\
&-\int_0^T ( \mu_f  (D \vel + {D \vel}^{\top}  ) \mnormalf - p \mnormalf + \mu_s  (D \mU + {D \mU}^{\top}  ) \mnormals \\ & \hspace{30pt} + \lambda_s  \mathrm{div}( \mU )\mnormals, \mpsiv  )_{\Gamma_i} d  \mt  = 0,
\end{align*}
which, by inserting the interface condition, yields the weak formulation
\begin{align*}
&\rho_f  \mdl \partial_\mt \vel, \mpsiv  \mdr_\mQf+ \int_0^T a_f (\vel, \mpsiv  ) d  \mt -  \mdl p, \mathrm{div}(\mpsiv)  \mdr_\mQf\\
&+\rho_s  \mdl \partial_\mt \vel, \mpsiv  \mdr_\mQs + \int_0^T a_s (\mU, \mpsiv ) d  \mt -  \mdl \mathrm{div}(\vel), \mpsip \mdr_\mQf\\
&-  \mdl \rho_f \mathbf f_f, \mpsiv \mdr_\mQf-  \mdl \rho_s \mathbf f_s, \mpsiv  \mdr_\mQs \\ &+ \mrhof ( \vel \mevalo - \vel_0, \mpsiv \mevalo  )_{\Omega_f}  + \mrhos ( \vel \mevalo - \vel_0, \mpsiv \mevalo  )_{\Omega_s} = 0.
\end{align*}
Hence, the linear PDE can be written as $\langle A(\metav, \metap), (\mpsiv, \mpsip) \rangle = \langle F, \mpsiv \rangle$, where the operator $A$ is defined by
\begin{align*}
\langle A ( \metav, \metap  ),  (\mpsiv, \mpsip ) \rangle &= \rho_f  \mdl \partial_\mt \metav, \mpsiv \mdr_{\mQf}+ \int_0^T a_f (\metav, \mpsiv  ) d  \mt -  \mdl \metap, \mathrm{div}(\mpsiv) \mdr_{\mQf}\\
& \hspace{-1.3cm} + \rho_s  \mdl \partial_\mt \metav, \mpsiv \mdr_{\mQs} + \int_0^T a_s ( \int_0^\mt \metav \mevals d \ms, \mpsiv  ) d  \mt -  \mdl \mathrm{div}(\metav), \mpsip \mdr_{\mQf} \\ & \hspace{-1.3cm} + \mrhof ( \metav \mevalo, \mpsiv \mevalo  )_{\Omega_f}  + \mrhos ( \metav \mevalo, \mpsiv \mevalo  )_{\Omega_s}.
\end{align*}
We are interested in the operator $A^*$.
The term which destroys the symmetry of the operator $A$ is given by $a_s$. A closer consideration of this term yields (under the assumption that we are on spaces where Fubini's theorem is valid):
\begin{align*}
& \int_0^T a_s (\int_0^\mt \metav  \mevals  d  \ms , \mpsiv  ) d  \mt = \int_0^T a_s (\int_0^\mt \metav  \mevals  d  \ms , \mpsiv  \mevalt  ) d  \mt \\
&= \int_0^T \int_0^\mt a_s (\metav \mevals, \mpsiv \mevalt )  d  \ms  d  \mt 
= \int_0^T \int_\ms^T a_s (\metav \mevals, \mpsiv \mevalt ) d  \mt  d  \ms  
\end{align*}
\begin{align*}
&= \int_0^T a_s (\metav \mevals, \int_\ms^T \mpsiv \mevalt d  \mt )  d  \ms  
= \int_0^T a_s (\metav \mevals, \int_0^{T-\ms} \mpsiv  (\cdot, T-\mt ) d  \mt  )  d  \ms  \\
&= \int_0^T a_s (\metav ( \cdot, T-\ms ), \int_0^\ms \mpsiv  (\cdot, T-\mt ) d  \mt  )  d  \ms 
\end{align*}
Introducing
$\mbetav (\xi, \mt) = \metav (\xi, T-\mt )$, $\mbetap(\xi, \mt) = \metap(\xi, T-\mt)$, $\mbpsiv (\xi, \mt) = \mpsiv (\xi, T-\mt )$, \newline $\mbpsip(\xi, \mt) = \mpsip(\xi, T-\mt)$ yields
\begin{align*}
\int_0^T a_s (\int_0^\mt \metav  \mevals d  \ms , \mpsiv  ) d  \mt = \int_0^T a_s (\int_0^\mt \mbpsiv  \mevals  d  \ms  ,\mbetav ) d  \mt.
\end{align*}
We have $\partial_\mt \mpsiv \mevalt = - \partial_\ms \mbpsiv (\cdot, \ms )$ for $\ms = T-\mt$ and the following equation holds true: 
\begin{align*}
&\mdl \metav, \partial_\mt \mpsiv  \mdr_{\mQf}= \int_0^T ( \metav  \mevalt, \partial_\mt \mpsiv  \mevalt  )_{\Omega_f} d  \mt \\
& = \int_0^T  ( \mbetav  (\cdot, T - \mt ), - \partial_\ms \mbpsiv  (\cdot,\ms)|_{\ms=T - \mt} )_{\Omega_f} d  \mt =  - \int_0^T  (  \mbetav  \mevals, \partial_\ms\mbpsiv \mevals )_{\Omega_f} d  \ms \\
& = - \int_0^T  ( \mbetav  \mevalt, \partial_\mt\mbpsiv \mevalt )_{\Omega_f} d  \mt 
= -  \mdl \mbetav, \partial_\mt \mbpsiv  \mdr_{\mQf}.
\end{align*}
Thus, integration by parts yields
\begin{align*}
\mdl \partial_\mt \metav, \mpsiv  \mdr_{\mQf} &= 
[( \metav (\cdot,\mt), \mpsiv (\cdot,\mt) )_{\Omega_f}]_{\mt=0}^{\mt=T}
-  \mdl \metav, \partial_\mt\mpsiv \mdr_{\mQf} \\
&=  ( \metav (\cdot, T ), \mpsiv (\cdot, T ) )_{\Omega_f} -  ( \metav \mevalo, \mpsiv \mevalo  )_{\Omega_f} +  \mdl \mbetav, \partial_\mt \mbpsiv  \mdr_{\mQf},\\
\mdl \partial_\mt \metav, \mpsiv  \mdr_{\mQs} &=  ( \metav (\cdot, T ), \mpsiv (\cdot, T )  )_{\Omega_s} - ( \metav \mevalo, \mpsiv \mevalo  )_{\Omega_s} +  \mdl \mbetav, \partial_\mt \mbpsiv \mdr_{\mQs}.
\end{align*}
Combining these results, rewriting the terms in $\mbetav$, $\mbetap$, $\mbpsiv$, $\mbpsip$ and using that \newline $( \metav(\cdot, T), \mpsiv(\cdot, T)  )_\Omega =  ( \mbetav\mevalo, \mbpsiv\mevalo  )_\Omega$ yields:
\begin{align*}
\langle A ( \metav, \metap  ),  (\mpsiv, \mpsip ) \rangle &= \rho_f  \mdl \partial_\mt \mbpsiv, \mbetav  \mdr_{\mQf} + \int_0^T a_f ( \mbpsiv, \mbetav  ) d  \mt -  \mdl \mathrm{div}(\mbpsiv), \mbetap \mdr_{\mQf} \\
&\hspace{-2.5cm}+ \rho_s  \mdl \partial_\mt \mbpsiv, \mbetav  \mdr_{\mQs} + \int_0^T a_s (\int_0^\mt \mbpsiv \mevals d  \ms  , \mbetav  ) d  \mt  -  \mdl \mbpsip, \mathrm{div}(\mbetav) \mdr_{\mQf} \\ & \hspace{-2.5cm} +  \mrhof (\mbpsiv\mevalo, \mbetav\mevalo  )_{\Omega_f} +  \mrhos (\mbpsiv\mevalo, \mbetav\mevalo  )_{\Omega_s} =  \langle A (\mbpsiv,\mbpsip ),  (\mbetav,\mbetap ) \rangle.
\end{align*}
Thus, the adjoint has the same structure as the forward model, but reverses the temporal flow of information.

\bibliographystyle{plain}
\bibliography{bibFSI}
\end{document}